\documentclass{amsart}
\usepackage{graphics,amssymb,amscd,amsthm,verbatim,psfrag}

\usepackage{graphics,verbatim}	

\newtheorem{proposition}{Proposition}[section]

\newtheorem{prop}[proposition]{Proposition}

\theoremstyle{definition}

\theoremstyle{remark}
\newtheorem{remark}[proposition]{Remark}

\newtheorem{nc_criterion}{Noncrossing Criterion}
\newtheorem{cl_criterion}{Compatibility Criterion}

\newcommand{\reals}{\mathbb R}

\newcommand{\set}[1]{{\left\lbrace #1 \right\rbrace}}

\newcommand{\ep}{\varepsilon}
\newcommand{\br}[1]{\langle #1 \rangle}
\newcommand{\Cat}{\mathbf{Cat}}
\newcommand{\Pge}{{\Phi_{\ge -1}}}

\newcommand{\figpage}[1]{
\special{#1}
\special{ps: /pageforpagesdotps exch def}
}

\begin{document}
\title[The Coxeter plane]{Noncrossing partitions, clusters and the Coxeter plane}

\author{Nathan Reading}
\address{
Department of Mathematics, North Carolina State University, Raleigh, North Carolina 27695-8205, USA}
\thanks{This project began while the author was partially supported by NSF grant DMS-0202430 and was completed while the author was partially supported by NSA grant H98230-09-1-0056.}
\email{nathan\_reading@ncsu.edu}
\urladdr{http://www4.ncsu.edu/$\sim$nreadin}

\begin{abstract}
When $W$ is a finite Coxeter group of classical type ($A$, $B$, or $D$), noncrossing partitions associated to $W$ and compatibility of almost positive roots in the associated root system are known to be modeled by certain planar diagrams.
We show how the classical-type constructions of planar diagrams arise uniformly from projections of small $W$-orbits to the \emph{Coxeter plane}.
When the construction is applied beyond the classical cases, simple criteria are apparent for noncrossing and for compatibility for $W$ of types $H_3$ and $I_2(m)$ and less simple criteria can be found for compatibility in types $E_6$, $F_4$ and $H_4$.
Our construction also explains why simple combinatorial models are elusive in the larger exceptional types.
\end{abstract}

\maketitle

\setcounter{tocdepth}{2}

\section{Introduction}\label{intro}

\subsection{$W$-Catalan combinatorics}
The ``classical'' noncrossing partitions are set partitions $\Pi$ of $[n+1]:=\set{1,2,\ldots n+1}$ such that when the elements of $[n+1]$ are placed in cyclic order on a circle, the convex hulls of the blocks of $\Pi$ are disjoint.  
Thus, for example, every partition of $[4]$ is noncrossing except $\set{\set{1,3},\set{2,4}}$, which is ``crossing'' because the line segment defined by $\set{1,3}$ intersects the line segment defined by $\set{2,4}$.
Noncrossing partitions were first studied by Kreweras~\cite{Kreweras} in 1972 and since then their rich combinatorial and enumerative structure has been widely studied. (Surveys include \cite[Chapter~4]{Armstrong} and~\cite{Simion}.)
Noncrossing partitions are counted by the Catalan number.
Recently, through the merger of lines of research in algebraic combinatorics~\cite{Ath-Rei,Reiner} and geometric group theory~\cite{Bessis,Biane,Bra-Wa} it became apparent that the classical noncrossing partitions are a special case (the case $W=S_{n+1}$) of a construction valid for any finite Coxeter group~$W.$
The $W$-noncrossing partitions are counted by the $W$-Catalan number $\Cat(W)$.

The general definition of noncrossing partitions is not phrased in terms of planar diagrams.
Rather, the parabolic subgroups of $W$ play the part of partitions, and an algebraic construction produces a set of parabolic subgroups which are to be considered ``noncrossing.''
When $W$ is the symmetric group, drawing the set partitions in the plane as described above, the noncrossing criterion of Kreweras and the algebraic criterion agree.
Similarly, when~$W$ is one of the other classical finite Coxeter groups ($B_n$ or $D_n$) there exist planar diagrams which can be used to determine whether a parabolic subgroup is crossing or noncrossing.

Another set of objects counted by the Catalan number has also recently been generalized to the context of finite Coxeter groups: triangulations of a convex polygon.
Fomin and Zelevinsky~\cite{ga} showed that a \emph{cluster algebra} satisfies a certain finiteness condition if and only if the underlying combinatorics of the cluster algebra is governed by a finite Coxeter group~$W$ (or more precisely, a finite crystallographic root system associated to~$W$).
The combinatorial structure of the cluster algebra is a simplicial complex whose vertex set is a subset (the \emph{almost positive roots}) of the root system.
The facets of the complex are \emph{clusters}: sets of almost positive roots that are pairwise ``compatible.''
The complex is dual to a simple polytope called the generalized associahedron.
The number of clusters, or equivalently the number of vertices of the generalized associahedron, is $\Cat(W)$.

When~$W$ is the symmetric group, the almost positive roots can be put in bijection with the diagonals of a convex polygon such that two almost positive roots are compatible if and only if the corresponding diagonals do not cross.
Thus a maximal set of compatible roots (or \emph{cluster}) corresponds to a maximal set of noncrossing diagonals of a polygon.
The latter is equivalent to a triangulation of the polygon.
When~$W$ is a Coxeter group of type $B_n$ or $D_n$, there are planar models which encode compatibility and which are only slightly more complicated.

Both the generalized associahedron and the noncrossing partition lattice have natural dihedral symmetry.
The planar models in types $A$, $B$ and $D$ realize these symmetries as dihedral symmetries acting geometrically on the plane.

\subsection{The main construction}\label{intro main}
This paper explains the planar diagrams in types $A$, $B$, and $D$ by giving a Coxeter-theoretically uniform construction of planar diagrams for parabolic subgroups and for almost-positive roots.
To the extent possible, we also give criteria for noncrossing and for compatibility (but not compatibility degree) in other types.
We also suggest a simple explanation for why planar diagrams are easy in types $A$, $B$, and $D$, but problematic in many exceptional types.

The definitions of clusters and noncrossing partitions involve the choice of a \emph{Coxeter element} for~$W$:
an element of~$W$ which can be expressed as the product of some permutation of the set $S$ of simple generators of~$W.$
The order of a Coxeter element is the \emph{Coxeter number} $h$ of~$W$ and the \emph{exponents} $e_1,\ldots,e_n$ of~$W$ are certain integers that can be read off from the eigenvalues of a Coxeter element.
(The Coxeter number and exponents are well-defined because any two Coxeter elements are conjugate in~$W.$)
The $W$-Catalan number is given by
\[\Cat(W)=\prod_{i=1}^n\frac{e_i+h+1}{e_i+1}.\]
The properties of Coxeter elements are closely tied to a certain ($2$-dimensional) plane $P$ which we call the \emph{Coxeter plane}, and which was first considered in generality by Coxeter~\cite{RegPoly}. 
A careful analysis of the Coxeter plane by Steinberg~\cite{Steinberg} provided the first uniform proofs of the key properties of Coxeter elements.

Our construction begins with orthogonal projections of a small $W$-orbit $o$ to the Coxeter plane.
Each parabolic subgroup $W'$ defines a partition of $o$ into $W'$-orbits.
The diagram for $W'$ is essentially this partition of $o$, projected to $P$.
The planar models for almost positive roots are obtained by altering the projected orbits, which have dihedral symmetry of order $2h$, to obtain a diagram with dihedral symmetry of order $2(h+2)$.

Once planar diagrams are constructed, it remains to find criteria for reading off, from the diagrams, whether a parabolic subgroup is noncrossing, or whether a pair of roots is compatible.
We give simple criteria for noncrossing and for compatibility for $W$ of types $H_3$ and $I_2(m)$, and slightly more complicated criteria for compatibility in types $E_6$, $F_4$ and $H_4$.
Although a uniform criterion for compatibility in all types remains elusive, the diagrams we construct for almost positive roots seem to come very close to having a uniform criterion.
Indeed, small \textit{ad hoc} alterations of the diagrams in the remaining types $E_6$, $E_7$, $E_8$ $F_4$, and $H_4$ produce a complete model for compatibility, with the simplest possible criterion.

The approach via the Coxeter plane indicates that the classical finite Coxeter groups (and now types $H_3$ and $I_2(m)$) admit simple criteria for noncrossing and for compatibility precisely because they admit small orbits.
In this context, a ``small'' orbit is an orbit whose size is approximately the Coxeter number $h$.
The groups of type $A_n$, $B_n$, $I_2(m)$, $D_n$, and $H_3$ each possess an orbit of size $h$ or $h+2$, while the remaining groups $F_4$, $E_6$, $E_7$, $H_4$, and $E_8$ have orbits of sizes $2h$, $2h+3$, $3h+2$, $4h$, and $8h$, respectively.
The intrinsic complexity of the diagrams unavoidably increases when the size of the orbit is much greater than $h$.
It is well-known that the property of having small orbits is also responsible for the combinatorial models which realize the classical groups as various groups of permutations, by letting the group act as permutations of a small orbit.  
However, in this latter context, a ``small'' orbit is one whose size is a small integer multiple of the rank of the Coxeter group.

\subsection{Computer-generated diagrams}\label{computer}
Postscript files are available on the author's website showing the results of the main construction for many finite Coxeter groups.
The figures were made with the help of John Stembridge's \texttt{coxeter} and \texttt{weyl} maple packages and PostScript drawing routines developed by the author.
The files are optimized for small file size as PostScript files, by taking advantage of the fact that PostScript is a complete programming language.
The author recommends viewing the files (or at least those files associated with the larger exceptional groups) with a PostScript viewer that interprets PostScript directly, rather than converting to PDF.
File sizes in PDF will be up to about 350 times larger than the postscript files.

Each file for parabolic subgroup diagrams (e.g.\ \texttt{e6nc.ps}) displays one representative from each symmetry class of diagram,  labels the parabolic subgroup as crossing or noncrossing, and gives a list of the reflections generating the parabolic subgroup.
The blocks in the partitions represented are indicated by line segments connecting pairs of points related by a reflection in the parabolic subgroup.
The blocks are also indicated by colors of vertices, with black points representing singleton blocks.
Colors of multiple points at the origin can be determined from the segments, keeping in mind that there are no degenerate segments connecting different points at the origin.
The reflections are numbered according to the numbering of positive roots in Stembridge's packages.
The noncrossing parabolic subgroups are identified by computer, via a bijection defined in~\cite[Section~6]{sortable} (see also \cite[Theorem~8.9]{typefree}) from \emph{sortable elements} to noncrossing parabolic subgroups.

Each file for diagrams of almost positive roots (e.g.\ \texttt{e6cl.ps}) pictures diagrams for pairs of almost positive roots $(\alpha,\beta)$, where $\alpha$ ranges over all negative simple roots and $\beta$ ranges over all almost positive roots.
Each diagram is a collection of line segments, labeled as necessary by which of the multiple origin points is an endpoint of the segment.
The labeling of segments is given in the diagram files by coloring the segments.
In Section~\ref{cl diag}, descriptions are given of \textit{ad hoc} alterations to the almost positive root diagrams in some types.  
These altered diagrams are also available, e.g.\ as \texttt{e6cl\_alt.ps}.

\subsection{Outline of the paper}
The remainder of the paper proceeds as follows.
In Section~\ref{orb sec}, we give background on crossing and noncrossing parabolic subgroups, relate these to partitions of a $W$-orbit, and describe the construction of Coxeter-planar diagrams for parabolic subgroups.
In Section~\ref{nc sec} we show that in the classical types, our construction of planar diagrams for partitions leads to the usual planar models for crossing/noncrossing partitions, and discuss how the construction works out in other types.
In Section~\ref{cl sec} we review the definition of compatibility and describe the usual models for compatibility in the classical types, in preparation for the development, in Section~\ref{root diag}, of a uniform construction of diagrams for almost positive roots.
Finally, in Section~\ref{cl diag}, we show that the construction recovers the classical planar diagrams, and discuss compatibility criteria in other types.

\section{Parabolic subgroups and partitions of orbits}\label{orb sec}
In this section we discuss crossing and noncrossing parabolic subgroups.
We introduce the notion of partitions of an orbit and describe how to obtain planar diagrams for parabolic subgroups by projecting partitions of an orbit to the Coxeter plane.

\subsection{Parabolic subgroups}
We assume basic background on Coxeter groups and root systems, including the classification of finite Coxeter groups.
This background can be found, for example, in~\cite{Bj-Br,Bourbaki,Humphreys}.
Throughout the paper, $(W,S)$ will stand for a finite Coxeter system of rank $n=|S|$.
For $s_1,s_2\in S$, the order of $s_1s_2$ will be denoted by $m(s_1,s_2)$.
We fix a representation of~$W$ as a real reflection group on a Euclidean space $V$, and make no distinction between 
elements of~$W$ and their action on $V$.
We also assume that the only element fixed by $W$ is the origin.
This implies in particular that the dimension of $V$ equals the rank of $W$.
The set of elements of $W$ which act as reflections is $T=\set{wsw^{-1}:s\in S, w\in W}$.
For each $t\in T$, let $H_t$ be the reflecting hyperplane associated to $t$.
The collection $\set{H_t:t\in T}$ is called the \emph{Coxeter arrangement} associated to $W$.

Throughout the paper, we will assume that $W$ is \emph{irreducible}, meaning that $W$ cannot be written as a direct product of Coxeter systems of strictly lower rank.
Equivalently, we require that the diagram of $W$ be connected as a graph.
The constructions of this paper rely, in two major ways, on the assumption that $W$ is irreducible.
First, the assumption of irreducibility, together with the assumption that $W$ fixes only the origin, implies that the linear span of any nontrivial orbit is $V$. 
Second, the assumption of irreducibility is essential to the construction of the Coxeter plane.

For any $J\subseteq S$, the subgroup $W_J$ generated by $J$ is called a \emph{standard parabolic subgroup}.
More generally, a \emph{parabolic subgroup} is any subgroup which is conjugate in $W$ to a standard parabolic subgroup.
The \emph{rank} of a parabolic subgroup $W'$ is $|J|$, where $W_J$ is some standard parabolic subgroup conjugate to $W'$. 
The parabolic subgroups of the symmetric group $S_{n+1}$ (the Coxeter group of type $A_n$) are exactly the subgroups~$W'$ generated by transpositions.
Such a subgroup~$W'$ decomposes the set $[n+1]$ into $W'$-orbits.
On the other hand, to each set partition $\Pi$ of $[n+1]$, one can associate the parabolic subgroup~$W'$ generated by transpositions $(i\,\,j)$ such that $i$ and $j$ are in the same block of $\Pi$.
Thus for an arbitrary finite Coxeter group $W,$ one may think of the parabolic subgroups of $W$ as a $W$-analog of partitions.

Another way to think of set partitions of $[n+1]$ is as subspaces in the intersection lattice of the Coxeter arrangement associated to $S_{n+1}$.
The Coxeter arrangement for the usual reflection representation of $S_{n+1}$ consists of all hyperplanes defined by equations $x_i=x_j$ for $1\le i<j\le n$.
Given a set partition $\Pi$ of $[n+1]$, one constructs a subspace by imposing the equation $x_i=x_j$ whenever $i$ and $j$ are in the same block of $\Pi$.
This is easily seen to be a bijection between set partitions and subspaces in the intersection lattice of $S_{n+1}$.

For general finite $W$, parabolic subgroups of $W$ are in bijection~\cite{Barcelo-Ihrig} with subspaces in the intersection lattice of the Coxeter arrangement for $W$, with a parabolic subgroup $W'$ mapping to the set stabilized by $W'$, a subspace whose codimension equals the rank of $W'$.
The inverse map takes a subspace in the intersection lattice to the subgroup of $W$ stabilizing that subspace.

\subsection{Crossing and noncrossing}\label{c nc sec}
We next explain the notion of crossing or noncrossing parabolic subgroups of $W,$ quoting results from~\cite{Bessis,Bra-Wa,BWlattice}.
Any element $w\in W$ can be written as a word in the alphabet $T$.
Since words in the alphabet $S$ are so commonly used in the study of Coxeter groups,  to avoid confusion we will refer a word in the alphabet $T$ as a \emph{$T$-word}.
A \emph{reduced} $T$-word for $w$ is a $T$-word for $w$ which has minimal length among all $T$-words for $w$.
The \emph{absolute length} of $w$ is the length of a reduced $T$-word for $w$.
(The usual \emph{length} of $w$ is the length of a reduced word for $w$ in the alphabet $S$.)
Define a partial order on~$W$ by setting $u\le v$ if and only if $u$ has a reduced $T$-word which is a prefix of some reduced $T$-word for $v$.
This partial order is called the \emph{absolute order} on $W$.

The interval $[1,c]$ in the absolute order, where $c$ is a Coxeter element, has come to be called the \emph{noncrossing partition lattice} associated to $W$.
For any element $x$ of $[1,c]$ and any reduced $T$-word $t_1\cdots t_i$ for $x$, the group $\br{t_1,\ldots,t_i}$ is a parabolic subgroup of~$W.$
Two different such words for $x\in [1,c]$ give rise to the same parabolic subgroup if and only if the products of the two words are the same element of~$W.$
The parabolic subgroups arising in the manner are called \emph{noncrossing parabolic subgroups}.
All other parabolic subgroups are \emph{crossing}.
The interval $[1,c]$ corresponds to containment order on noncrossing parabolic subgroups or, equivalently, to reverse containment order on the subspaces associated to noncrossing parabolic subgroups.
If~$W$ is reducible as a direct product of nontrivial Coxeter groups, then the noncrossing parabolic subgroups of~$W$ are exactly those parabolic subgroups whose intersection with each factor is a noncrossing parabolic subgroup in the factor.
Thus it is harmless (and essential, as discussed above) to restrict to the case where~$W$ is irreducible.

\begin{remark}\label{terminology}
Although $[1,c]$ is called the noncrossing partition lattice, it makes sense to reserve the term ``noncrossing partitions'' for noncrossing parabolic subgroups:
The elements of $[1,c]$ are group elements, not partitions in any natural sense, and there are no analogous group elements playing the role of ``crossing partitions.''
\end{remark}

Conjugation by an element of $W$ permutes $T$ and thus acts as an isomorphism of the absolute order.
Since any two Coxeter elements are conjugate in $W$, the isomorphism type of $[1,c]$ does not depend on the choice of $c$.
The set of elements in $[1,c]$ and the set of noncrossing parabolic subgroups do depend on the choice of $c$, but changing from $c$ to another Coxeter element $wcw^{-1}$ affects both of these sets by a global conjugation by $w$.
It will be useful to fix a particularly convenient Coxeter element.
The Coxeter diagram of a finite Coxeter group is a tree, and in particular is bipartite.
Let $S_+\cup S_-$ be a bipartition of the Coxeter diagram for~$W.$
Define involutions 
\[c_+:=\prod_{s\in S_+}s\ \ \mbox{ and }\ \ c_-:=\prod_{s\in S_-}s\]
so that $c:=c_-c_+$ is a Coxeter element.
A Coxeter element arising in this manner is called a \emph{bipartite Coxeter element}.
Since the elements of $S_+$ commute pairwise, the product $c_+$ is well-defined, and similarly $c_-$ is well-defined.
All further references to noncrossing partitions or crossing/noncrossing parabolic subgroups will refer to this fixed bipartite Coxeter element.
The noncrossing parabolic subgroups for $c$ and $c^{-1}$ coincide, and thus since $c^{-1}=c_-cc_-=c_+cc_+$, the conjugation action of the subgroup $\br{c_+,c_-}$ preserves the crossing/noncrossing status of parabolic subgroups.

\subsection{Partitions of orbits}
Besides parabolic subgroups and subspaces in the intersection lattice, another $W$-analog of set partitions of $[n+1]$ appears not to have been considered before.
Given a fixed non-trivial orbit $o$ of $W$, each parabolic subgroup~$W'$ decomposes $o$ into $W'$-orbits.
We will call this decomposition the partition of $o$ induced by $W'$.
Collectively, such partitions of $o$ will be called $W$-partitions of $o$.
For general finite $W$, given a partition of $o$ that is known to the be a $W'$-partition for some~$W'$, one can read off~$W'$ from the partition of $o$.
Specifically, it is the stabilizer of its fixed subspace, and the latter is read off as in the following proposition, whose simple proof is omitted.
\begin{prop}\label{read off}
Let $o$ be a $W$-orbit and let $\Pi$ be a partition of $o$ induced by some parabolic subgroup $W'$ of $W$.
Then the fixed subspace of $W'$ is the orthogonal complement of the linear span of all vectors $x-y$ where $x,y\in o$ are in the same block of the partition.
\end{prop}

In light of Proposition~\ref{read off}, the $W$-partitions of $o$ may be thought of as a $W$-analog of the usual partitions of $[n+1]$.
Note that, even fixing $W=S_{n+1}$, the $W$-partitions of orbits generalize the usual partitions of $[n+1]$.
The usual partitions of $[n+1]$ arise when $o$ is taken to be a smallest nontrivial $S_{n+1}$-orbit.

\subsection{Diagrams for partitions of orbits} \label{par sec}
Our approach to combinatorial models of noncrossing partitions is to fix a $W$-orbit $o$, associate a planar diagram to each $W$-partition of $o$, and then seek a combinatorial characterization of crossing and noncrossing partitions of $o$ in terms of the planar diagrams.

The planar diagrams will arise from the orthogonal projection of $o$ to the \emph{Coxeter plane} $P$.
Details on~$P$ can be found in \cite{RegPoly}, in~\cite{Steinberg}, in Sections 3.16--3.20 of~\cite{Humphreys} or in Section V.6.2 of~\cite{Bourbaki}.
The important properties of~$P$ are as follows:
$P$ is fixed as a set by the action of $\br{c_+,c_-}$, and $\br{c_+,c_-}$ acts on~$P$ as a dihedral reflection group of order~$2h$.
Taking $\ep\in\set{+,-}$, the transformation $c_\ep$ stabilizes $H_s$ for $s\in S_\ep$.
In $P$ there are two lines $L_+$ and $L_-$ with the following properties:
A reflecting hyperplane $H$ for a reflection in $W$ intersects~$P$ in $L_\ep$ if and only if $H$ is the reflecting hyperplane for some simple reflection $s\in S_\ep$.
The other reflecting hyperplanes intersect~$P$ in lines which are the images of $L_\pm$ under $\br{c_+,c_-}$.
In particular, every reflecting hyperplane is in the $c$-orbit of $H_s$ for some $s\in S$.
Equivalently, every reflection in~$W$ is in the orbit (under conjugation by $c$) of some simple reflection $s$ and, more specifically, we have the following proposition.
(Cf. \cite[Proposition~V.6.2]{Bourbaki} or \cite[Exercise~3.19.2]{Humphreys}).

\begin{prop}
\label{c orbits}
For any $t\in T$, the orbit of $H_t$ under the action of the dihedral group $\br{c_+,c_-}$ either:
\begin{enumerate}
\item[(i) ]has $h/2$ elements and intersects $\set{H_s:s\in S}$ in a single element.
\item[(ii) ]has $h$ elements and intersects $\set{H_s:s\in S}$ in a two-element set.
\end{enumerate}
Furthermore, when (i) holds, the hyperplane in the intersection has $H=w_0H$, where $w_0$ is the longest element of $W$.  
When (ii) holds, the two elements in the intersection are related by the action of $w_0$.
\end{prop}

Fix a $W$-orbit $o$, and consider the orthogonal projection of $o$ to $P$.
The projected orbit is thought of as a point configuration, where some points have \emph{multiplicity} greater than one because the orthogonal projection to~$P$ is not always one-to-one.
Multiple points are considered to be distinguishable---each copy of a point is formally labeled by its pre-image in $o$.
The problem of multiplicity can be minimized by taking $o$ as small as possible, in which case the only point of the projected orbit having multiplicity is the origin, with multiplicity at most 3.  

\begin{figure}
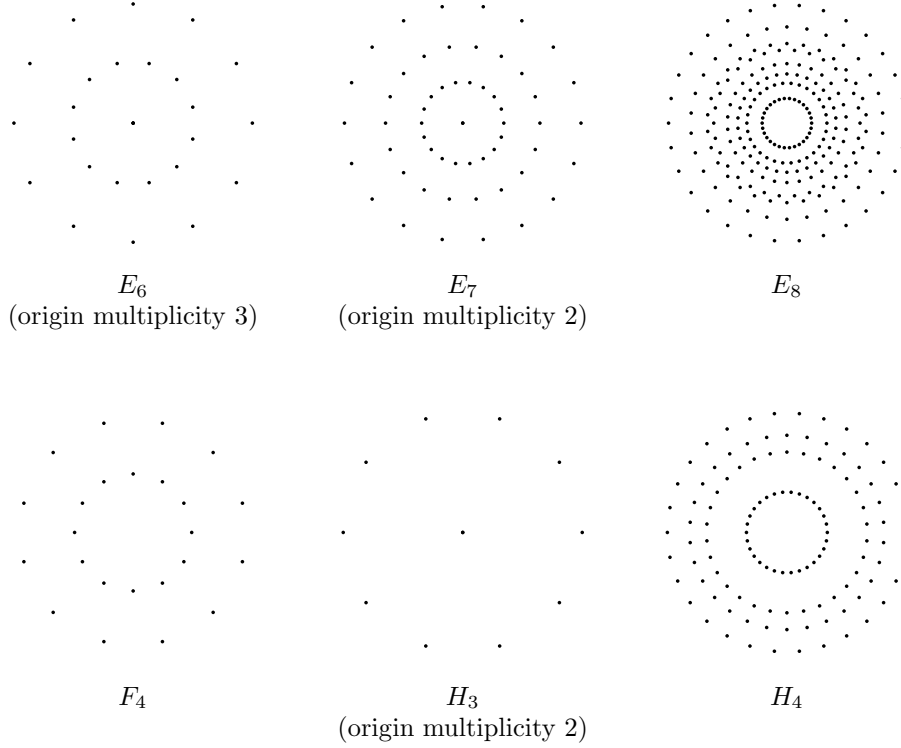

\begin{center}
\begin{tabular}{ccccccc}
\scalebox{0.9}{\includegraphics{e6orb.ps}}
&&&\scalebox{0.9}{\includegraphics{e7orb.ps}}
&&&\scalebox{0.9}{\includegraphics{e8orb.ps}}\\[8 pt]
$E_6$&&&$E_7$&&&$E_8$\\
(origin multiplicity 3)&&&(origin multiplicity 2)
\\[30 pt]
\scalebox{0.9}{\includegraphics{f4orb.ps}}
&&&\scalebox{0.9}{\includegraphics{h3orb.ps}}
&&&\scalebox{0.9}{\includegraphics{h4orb.ps}}\\[8 pt]
$F_4$&&&$H_3$&&&$H_4$\\
&&&(origin multiplicity 2)
\end{tabular}
\end{center}
\caption{The smallest orbit of $W,$ projected orthogonally to the Coxeter plane.}
\label{proj orb}
\end{figure}

For~$W$ of type $A$, $B$ or $I_2$, the smallest orbit is of size $h$, so the projection onto~$P$ consists of the vertices of an $h$-gon.
The smallest orbit in type $D_n$ or $H_3$ is of size $h+2$, and the projection onto~$P$ consists of an $h$-gon and two points at the origin.
These assertions about projected orbits in types $A$, $B$, and $D$ are argued in Sections~\ref{nc A}, \ref{nc B} and~\ref{nc D}.
For the exceptional types, the projected orbits are determined computationally.
Figure~\ref{proj orb} shows, for each exceptional Coxeter group $W,$ the projection of a smallest orbit of~$W$ onto~$P$.
In the cases where points in the projected orbit project to the origin, the multiplicity of points at the origin is indicated in the figure.
No other points in the projected orbits have multiplicity.

Figure~\ref{proj orb} makes it clear why, as mentioned in Section~\ref{intro main}, a ``small'' orbit is an orbit whose size is on the order of the Coxeter number $h$:
Since the Coxeter element $c$ is in particular an orthogonal transformation and fixes~$P$ as a set, it commutes with orthogonal projection to~$P$, so that the projected orbit decomposes as the union of $c$-orbits in~$P$.
When the size of the $W$-orbit is on the order of $h$, the projected $W$-orbit consists of a single $c$-orbit (forming a regular $h$-gon in~$P$) and perhaps several points at the origin.

Once the projection of $o$ to $P$ is calculated, the \emph{Coxeter-plane diagram} for each parabolic subgroup~$W'$ is obtained by projecting the partition of $o$ induced by $W'$ orthogonally to $P$.
The diagram can be thought of as a coloring of the projection of $o$ to~$P$, where the colors are the blocks of the partition of $o$.
Repeated points at the origin may have different colors.

\section{Diagrams for noncrossing partitions} \label{nc sec}
In this section, we discuss, for $W$ of various types, the construction of Coxeter-plane diagrams for parabolic subgroups.
Computer data are available for many types, as described in Section~\ref{computer}.

\subsection{Type $A_n$}\label{nc A}
We represent the Coxeter group $W$ of type $A_n$ in the usual way as the symmetric group $S_{n+1}$ acting on $\reals^{n+1}$ by permuting coordinates.
To obtain a representation fixing only the origin, one takes the quotient of $\reals^{n+1}$ modulo the line spanned by the vector $(1,1,\ldots,1)$.
We will describe points in $\reals^{n+1}$ and tacitly assume the quotient.
A smallest orbit $o$ of $W$ is the set of unit basis vectors $\set{e_1,\ldots,e_{n+1}}$.
The simple generators are $\set{s_i:i\in[n]}$, with $s_i$ acting by transposing $e_i$ and $e_{i+1}$.
Choose the bipartition of the diagram with $S_+=\set{s_1,s_3,\ldots}$ and $S_-=\set{s_2,s_4,\ldots}$, so that the bipartite Coxeter element $c=c_-c_+$ is the $(n+1)$-cycle $c=(1,3,\ldots,n+1,n,n-2,\ldots,2)$ if $n$ is even or $c=(1,3,\ldots,n,n+1,n-1,\ldots,2)$ if $n$ is odd.

The orbit $o$ is also a single orbit under the action of $c$.
This implies in particular that no point in $o$ projects orthogonally to the origin in $P$:
If some $e_i$ projects to the origin, then since $c$ fixes $P$ and thus fixes $P^\perp$, each element of $o$ is contained in $P^\perp$.
This is forbidden by considerations of dimension.
We conclude that the projection of $o$ to $P$ is a regular $(n+1)$-gon (i.e.\ an $h$-gon) centered at the origin.
The cyclic order on the projections of the vertices corresponds to the $(n+1)$-cycle $c$.
The $W$-partitions of $o$ coincide with the set partitions of $o$.
It is well-known that the noncrossing partitions of type $A$ are exactly the noncrossing partitions of a cycle introduced by Kreweras~\cite{Kreweras}. 
The proof can be found in~\cite{Biane} or~\cite{Brady}.
Thus noncrossing parabolic subgroups are identified by the following criterion.
\begin{nc_criterion} \label{A}
A parabolic subgroup is noncrossing if and only if the convex hulls of the blocks of its Coxeter-plane diagram are disjoint.
\end{nc_criterion}
To correctly compare the Coxeter-plane construction with the arguments in~\cite{Biane} or~\cite{Brady}, one should construct the noncrossing partitions of the cycle $c$ described above, rather than the usual noncrossing partitions of $(1,2,\ldots,n+1)$.
The latter $(n+1)$-cycle is a different Coxeter element $s_1s_2\cdots s_n$, which is not bipartite with respect to our choice of simple system.

\subsection{Type $B_n$}\label{nc B}
The Coxeter group of type $B_n$ is realized as the symmetry group of the $n$-dimensional cube with vertices $(\pm1,\pm1,\ldots,\pm1)$.
The set $o=\set{\pm e_1,\ldots,\pm e_n}$ is a smallest orbit.
The simple generators are $\set{s_i:i\in[n]}$, where $s_n$ is the reflection orthogonal to $e_n$, and $s_i$ acts by transposing $e_i$ and $e_{i+1}$, for each $i\in[n-1]$.
Choose $S_+=\set{s_1,s_3,\ldots}$ and $S_-=\set{s_2,s_4,\ldots}$, and define the bipartite Coxeter element $c=c_-c_+$.
When $n$ is even, 
\[c=(1,3,\ldots,n-1,-n,-n+2,\ldots,-2,-1,-3,\ldots,-n+1,n,n-2,\ldots,2),\]
and when $n$ is odd,
\[c=(1,3,\ldots,n,-n+1,-n+3,\ldots,-2,-1,-3,\ldots,-n,n-1,n-3,\ldots,2).\]
Similarly to the type-$A$ case, no element of $o$ projects to the origin of $P$, and since the antipodal pairs in $o$ must project to antipodal pairs in $P$, we conclude that the projection of $o$ to $P$ is a regular $(h=2n)$-gon centered at the origin.
Coxeter-plane diagrams of parabolic subgroups are centrally symmetric set partitions of the vertices of the $2n$-gon.
Noncrossing Criterion~\ref{A} is the combinatorial crossing/noncrossing criterion of~\cite{Reiner}, which has been shown in \cite{Bessis,BGN,Bra-Wa} to correctly classify parabolic subgroups as crossing/noncrossing in the algebraic sense.
Again, to correctly compare the above construction to \cite{Bessis,BGN,Bra-Wa,Reiner}, one must take the correct $2n$-cycle $c$ described above.

\subsection{Type $I_2(m)$}\label{nc I}
For a dihedral Coxeter group $W$ of type $I_2(m)$, the parabolic subgroups are $\set{1,t}$ for each reflection $t$ as well as the trivial subgroup $\set{1}$ and the entire group $W.$
In this case, $h=m$ and the plane~$P$ is the plane on which $W$ acts (so that projection to~$P$ is the identity map).
The orbit $o$ can be taken to be the vertices of a regular $m$-gon.
For any finite Coxeter group $W,$ the trivial subgroup, the subgroups $\set{1,t}$ for reflections $T$ and the whole group $W$ are noncrossing parabolic subgroups.
Thus in type $I_2(m)$ we must declare every parabolic subgroup to be noncrossing.
This can be accomplished with Noncrossing Criterion~\ref{A}.
The validity of the criterion is trivial for $W'=\set{1}$ and $W'=W$.
In the remaining cases, the criterion is valid because any parabolic subgroup $\set{1,t}$ decomposes $o$ into a collection of non-overlapping parallel line segments. 

\subsection{Type $D_n$}\label{nc D}
The Coxeter group of type $D_n$ can be realized in $\reals^n$ as the Coxeter group generated by simple reflections $s_i$ transposing $e_i$ and $e_{i+1}$ for $i\in[n-1]$ and $s_n$ transposing $e_{n-1}$ and $-e_n$.
A smallest orbit is $o=\set{\pm e_1,\ldots,\pm e_n}$.
We bipartition the simple reflections, except for $s_n$, by $S_+=\set{s_1,s_3,\ldots}$ and $S_-=\set{s_2,s_4,\ldots}$, and put $s_n$ is the part with $s_{n-1}$.
The resulting bipartite Coxeter element $c=c_-c_+$ for $n$ even is
\[(1,3,\ldots,n-1,-n+2,-n+4,\ldots,-2,-1,-3,\ldots,-n+1,n-2,n-4,\ldots,2)\,(n,-n),\]
and for $n$ odd is
\[(1,3,\ldots,n-2,-n+1,-n+3,\ldots,-2,-1,-3,\ldots,-n+2,n-1,n-3,\ldots,2)\,(n,-n).\]
If an element $x$ of $o$ does not project to the origin of $P$, then the projection of $x$ must have an $h$-fold orbit under the action of $c$.
(Here, the Coxeter number $h$ is $2n-2$.)
However, since orthogonal projection to $P$ commutes with the action of $c$, the point $x$ also has an $h$-fold orbit.
Since $\set{\pm e_n}$ is a $c$-orbit, we conclude that the elements $\pm e_n$ both project to the origin.
The other $2n-2$ elements of $o$ project to the vertices of a regular $h$-gon in $P$.

In this example, we catch the first glimpse of the complications inherent in the construction:
Since several points of the projected orbit coincide, it is unavoidable to have some intersections of blocks, even for noncrossing partitions.
However, the complications are tame in this case, and the following criterion correctly classifies parabolic subgroups as crossing or noncrossing.
\begin{nc_criterion} \label{D}
Two blocks in a Coxeter-plane diagram are \emph{noncrossing} if their convex hulls either are disjoint or intersect in a single point on the boundary of both.
A parabolic subgroup of $o$ is noncrossing if and only if the blocks in its Coxeter-plane diagram are pairwise noncrossing.
\end{nc_criterion}
\noindent
Noncrossing Criterion \ref{D} is a rephrasing of the criterion given by Athanasiadis and Reiner in~\cite[Section~3]{Ath-Rei}, again taking into account the labeling of the $(2n-2)$-cycle given by $c$ above.
Although Noncrossing Criterion~\ref{D} has a formally weaker requirement for the noncrossing property, it is still valid in types $A$, $B$, and $I$, because in those types any two blocks that intersect necessarily intersect in their relative interiors.

\subsection{Type $H_3$}\label{nc H3}
\begin{figure}
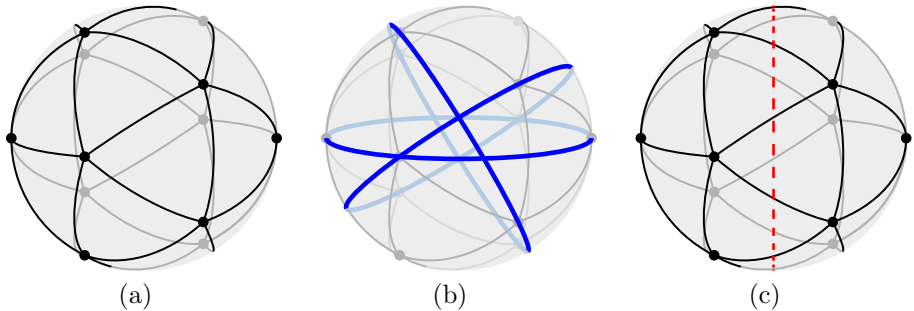

\begin{tabular}{ccc}
\figpage{ps: 1}\includegraphics{h3planes.ps}&
\figpage{ps: 2}\includegraphics{h3planes.ps}&
\figpage{ps: 3}\includegraphics{h3planes.ps}\\
(a)&(b)&(c)
\end{tabular}
\caption{Simple reflections and Coxeter plane for $H_3$.}
\label{h3icos}
\end{figure}

The Coxeter group $W$ of type $H_3$ is the symmetry group of the icosahedron and of the dodecahedron.
The smallest nontrivial orbit $o$ of $W$ consists of the 12 vertices of the icosahedron.
We will relate all of the $H_3$ constructions to the icosahedron depicted in Figure~\ref{h3icos}.a as a triangulation of a round sphere.
The sphere is partially transparent, so that we see the vertices and edges on the ``back'' of the sphere in gray.
Figure~\ref{h3icos}.b shows a choice of simple reflections for $W$, with each reflection represented by the intersection of its reflecting plane with the sphere.
(The vertices and edges of the icosahedron are faded in Figure~\ref{h3icos}.b to let the reflections stand out.)
The reflection whose plane is roughly horizontal in the picture will be called $s_2$, while $s_1$ and $s_3$ are the reflections with $m(s_1,s_2)=3$ and $m(s_2,s_3)=5$.
Choose $S_+=\set{s_1,s_3}$ and $S_-=\set{s_2}$, so that $c=c_-c_+=s_2s_1s_3$ is a bipartite Coxeter element.
The Coxeter plane~$P$ is indicated by the dotted line in Figure~\ref{h3icos}.c, and the Coxeter element $c$ acts by reflecting through~$P$ and rotating~$P$ by $1/10$ of a turn.
(The rotation moves points on the ``front'' of the sphere downwards in the picture.)

Two opposite vertices of the icosahedron project to the origin in~$P$.
The other 10 points project to the vertices of a regular decagon centered at the origin.
The noncrossing property of parabolic subgroups is described by Noncrossing Criterion~\ref{D}.

\begin{figure}
\centerline{\includegraphics{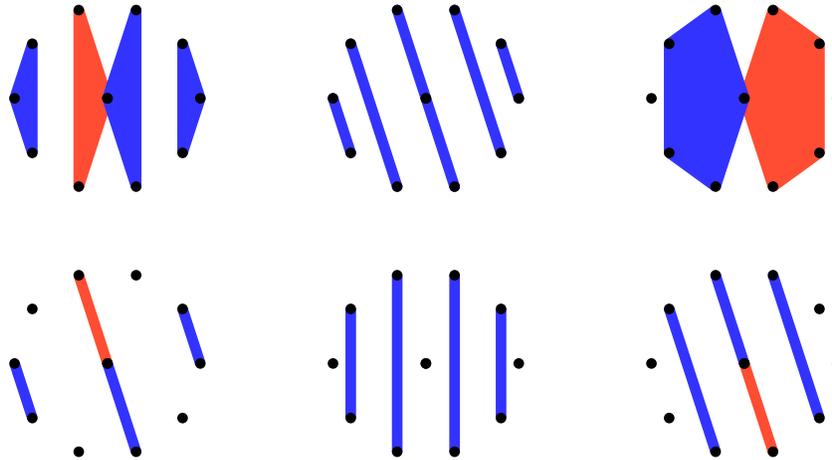}}
\caption{Noncrossing partitions in $H_3$.}
\label{h3nc}
\end{figure}

Figure~\ref{h3nc} shows the Coxeter-plane diagram for one representative of each $c$-conjugation orbit of noncrossing parabolic subgroups, omitting $\set{1}$ and~$W.$
Call the points at the origin $\pm 0$.
If a diagram has a block shown in red (or the lighter of the two grays) then the red block contains $-0$ but not $+0$.
The $c$-conjugation orbit of each noncrossing parabolic subgroup shown is of cardinality 5.
To see that the cardinality is 5 and not 10 in some cases, recall that the action of $c$ on the diagrams not only rotates by $1/10$ of a turn, but also switches the two points at the origin.
The bottom row of Figure~\ref{h3nc} shows noncrossing diagrams for parabolic subgroups of the form $\set{1,s}$ for a simple reflection $s$.
Reading from left to right in this bottom row, the simple reflections are $s_1$, $s_2$ and $s_3$.

\begin{figure}
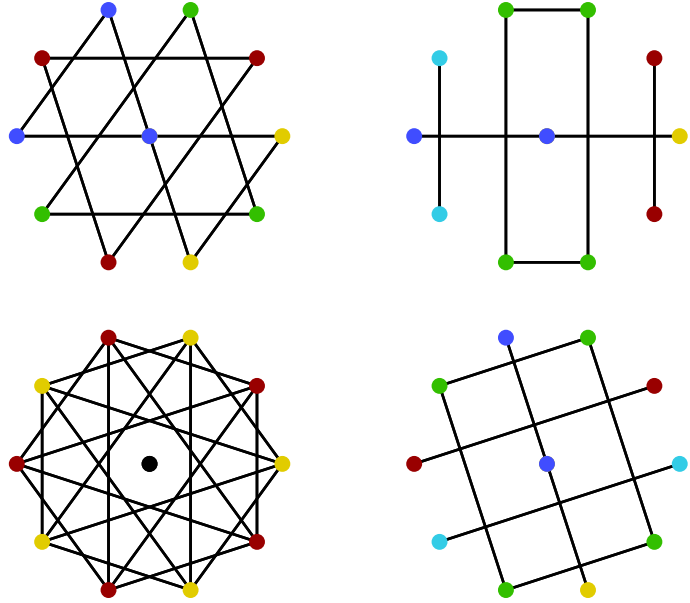

\begin{tabular}{ccccc}
\scalebox{1}{\figpage{ps: 2}\includegraphics{h3c.ps}}&&&&
\scalebox{1}{\figpage{ps: 1}\includegraphics{h3c.ps}}\\[20 pt]
\scalebox{1}{\figpage{ps: 4}\includegraphics{h3c.ps}}&&&&
\scalebox{1}{\figpage{ps: 3}\includegraphics{h3c.ps}}
\end{tabular}
\caption{Crossing parabolic subgroups of $H_3$.}
\label{h3c}
\end{figure}

Figure~\ref{h3c} shows one representative of each $c$-conjugation orbit of crossing parabolic subgroups.
Instead of the convex hulls of the blocks, the figure shows line segments indicating that two points are related by a reflection in the parabolic subgroup.
The partition of $o$ induced by $W'$ is thus the closure of the relation shown by the line segments.
The blocks are also indicated by colors of vertices.
The bottom-left diagram has two singleton blocks at the origin (shown in black).
The other three diagrams have both a blue and a yellow vertex at the origin.

The standard parabolic subgroup $W_{s_1,s_2}$ decomposes $o$ into four orbits of three elements each, arranged into four equilateral triangles on four parallel planes.
Two of the equilateral triangles are faces of the icosahedron.
Of the ten parabolic subgroups conjugate to $W_{s_1,s_2}$, five are noncrossing and five are crossing.
One of the noncrossing parabolic subgroups is $W_{s_1,s_2}$ itself, which appears as the top-left diagram in Figure~\ref{h3nc}.
A crossing parabolic subgroup conjugate to $W_{s_1,s_2}$ is shown as the top-left diagram in Figure~\ref{h3c}.

The standard parabolic subgroup $W_{s_1,s_3}$ decomposes $o$ into four orbits of two elements each, and one orbit containing 4 elements.
The four-element orbit consists of the vertices of two antipodal edges of the icosahedron.
Of the fifteen parabolic subgroups conjugate to $W_{s_1,s_3}$, five are noncrossing and ten are crossing, with $W_{s_1,s_3}$ appearing as the top-center diagram in Figure~\ref{h3nc}, and crossing parabolic subgroups conjugate to $W_{s_1,s_3}$ appearing as the two diagrams on the right in Figure~\ref{h3c}.

The standard parabolic subgroup $W_{s_2,s_3}$ decomposes $o$ into two orbits of five elements each, arranged into regular pentagons on two parallel planes, and two singleton orbits.
Of the six parabolic subgroups conjugate to $W_{s_1,s_2}$, five are noncrossing (the top-right diagram in Figure~\ref{h3nc} is $W_{s_1,s_2}$) and one is crossing (the bottom-left diagram in Figure~\ref{h3c}).

\subsection{Type $F_4$}\label{nc F}
\begin{figure}
\centerline{\scalebox{1}{\includegraphics{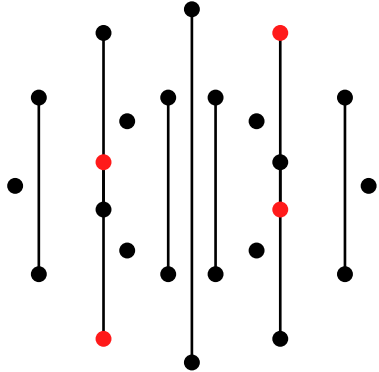}}}
\caption{A (noncrossing) parabolic subgroup of rank 1 in $F_4$.}
\label{F4 bad ref}
\end{figure}

The Coxeter group $W$ of type $F_4$ is the smallest example of an irreducible Coxeter group having no orbit of cardinality approximately $h$.
A smallest orbit $o$ of $F_4$ has $24$ elements, while $h$ is $12$.
(This orbit consists of the vertices of a regular $24$-cell.)
The projection of $o$ to~$P$ consists of two regular $h$-gons centered at the origin as shown in Figure~\ref{proj orb}.
The presence of two $h$-gons in the Coxeter-plane diagrams for parabolic subgroups causes problems for Noncrossing Criteria~\ref{A} and~\ref{D} which appear already in the case where the parabolic subgroup is generated by a single reflection.
Any parabolic subgroup generated by a single reflection is noncrossing.
However, Figure~\ref{F4 bad ref} shows the diagram of such a parabolic subgroup which fails both Noncrossing Criterion~\ref{A} and Noncrossing Criterion~\ref{D}.
The decomposition into blocks is finer than the decomposition of the diagram into connected pieces.
Specifically, two of the connected pieces of the diagram are further divided into blocks of two points each, as indicated by colors of points.
Thus in two places in the diagram, convex hulls of blocks intersect in their relative interiors.

\begin{figure}
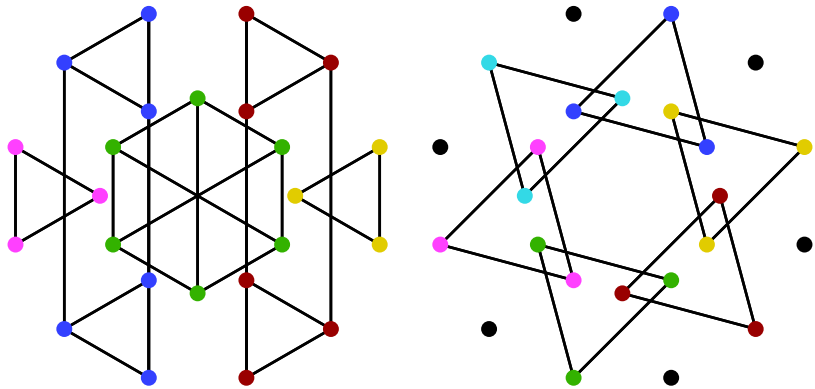

\begin{tabular}{cc}
\scalebox{1}{\includegraphics{f4cr.ps}}&\scalebox{1}{\includegraphics{f4non.ps}}
\end{tabular}
\caption{a: A crossing parabolic subgroup in $F_4$.  b:  A noncrossing parabolic subgroup in $F_4$.}
\label{F4 diff}
\end{figure}

To further illustrate the difficulty of finding a noncrossing criterion for Coxeter plane diagrams of type $F_4$, observe for example that the criterion must declare Figure~\ref{F4 diff}.a to be crossing while declaring Figure~\ref{F4 diff}.b to be noncrossing.
In Figure~\ref{F4 diff}, as in Figure~\ref{h3c}, line segments indicate that two points are related by a reflection in the parabolic subgroup and the blocks are also indicated by colors of vertices, with black dots indicating singleton blocks.
Various unsatisfying criteria can be devised to detect crossing parabolic subgroups in~$F_4$, for example in terms of centers of mass of projected blocks.
However, even these contrived conditions fail dramatically for small parabolic subgroups in $H_4$.

\section{Clusters} \label{cl sec}
In this section, we give a brief account of \emph{compatibility} and \emph{clusters of almost positive roots}.
These were introduced by Fomin and Zelevinsky to describe the vertices of the generalized associahedron~\cite{ga} and also to describe compatibility of \emph{denominator vectors} of \emph{cluster variables} in a \emph{cluster algebra}~\cite{ca2}.
The results reviewed in this section are all contained in~\cite{ga}.
We continue under the vital assumption that $W$ is irreducible.

Let $\Phi$ be a root system for $W$, chosen to be consistent with the fixed reflection representation.
Write $\Pi$ for the simple roots (corresponding to the simple generators $S$).
Let $\Phi^+$ denote the positive roots, which are in one-to-one correspondence with the reflections $T$.
Let $\Pge=\Phi_+\cup(-\Pi)$ be the set of \emph{almost positive roots} in the root system $\Phi$ associated to $W$.
As before, let $S=S_+\cup S_-$ be a bipartition of the Coxeter diagram, and define involutions $\tau_+$ and $\tau_-$ on $\Pge$:
\[\tau_\varepsilon(\alpha) = 
\begin{cases} 
\ \ \alpha & \mbox{if }\alpha=-\alpha_s\mbox{ for }s \in S_{-\ep}\\
c_\ep(\alpha) & \text{otherwise.} 
\end{cases} 
\end{equation*} 

The following is a direct consequence\footnote{Proposition~\ref{tau orbits} was first pointed out as~\cite[Theorem~2.6]{ga}, where its proof relies on \cite[Lemma~2.1]{ga}=\cite[Exercise~V.6.2]{Bourbaki}, which relies on the considerations leading to Proposition~\ref{c orbits}.} of Proposition~\ref{c orbits}.

\begin{prop}
\label{tau orbits}
For any $\alpha\in\Pge$, the orbit of $\alpha$ under the action of the dihedral group $\br{\tau_+,\tau_-}$ either:
\begin{enumerate}
\item[(i) ]has $(h+2)/2$ elements and intersects $-\Pi$ in a single element.
\item[(ii) ]has $h+2$ elements and intersects $-\Pi$ in a two-element set.
\end{enumerate}
When (i) holds, the root in the intersection has $\beta=-w_0\beta$, where $w_0$ is the longest element of $W$.  
When (ii) holds, the two roots in the intersection are related by the action of $-w_0$.
\end{prop}

There is a unique binary relation (called ``compatibility'') on $\Phi_{\geq -1}$ which is both $\langle \tau_-,\tau_+ \rangle$-invariant and has the following property for each negative simple root $-\alpha_s$:
A root $\beta\in\Pge$ is compatible with $-\alpha_s$ if and only if $\alpha_s$ appears with coefficient zero in the simple root expansion of $\beta$.
(More generally, the paper~\cite{ga} defines the the \emph{compatibility degree} $(\alpha\parallel\beta)$ of a pair $\alpha,\beta$ of almost positive roots, which we will not consider here.)
The compatibility relation is symmetric.

Maximal sets of pairwise-compatible almost positive roots are called \emph{clusters}.
Each cluster contains $n$ elements, where $n$ is the rank of $W$.
The clusters are the maximal faces of a simplicial complex called the \emph{cluster complex}, with faces given by pairwise compatible sets of almost positive roots.
Dual to the cluster complex is a simple polytope called the \emph{generalized associahedron} for $W$.

In Section~\ref{root diag}, we describe a construction of a planar diagram for each almost positive root, with the goal of describing compatibility of two roots in terms of the union of their two diagrams.
For the purposes of comparison and motivation, we conclude this section by describing the usual planar diagrams for compatibility in types $A_n$, $B_n$ and $D_n$ as described in \cite[Section~3.5]{ga}.

\subsection{The usual model for compatibility in type $A_n$}\label{usual A}
Let $W$ be a Coxeter group of type $A_n$ with the bipartition of $S$ described in Section~\ref{nc A}.
For each $i\in[n]$, let $\alpha_i$ be the simple root corresponding to $s_i$.

\begin{figure}
\begin{center} 
\scalebox{.95}{
\begin{picture}(0,0)(0,0) 
\put(26,60){$-\alpha_1$}
\put(36,130){$-\alpha_2$}
\put(78,60){$-\alpha_3$}
\put(95,130){$-\alpha_4$}
\put(140,65){$-\alpha_5$}
\put(155,140){$-\alpha_6$}
\put(187,95){--$\alpha_7$}
\end{picture} 
\includegraphics{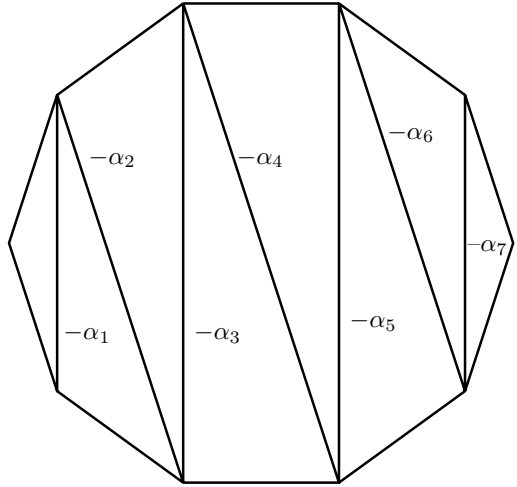}
}
\end{center} 
\caption{The ``snake'' in type $A_7$} 
\label{Asnake} 
\end{figure}

Compatibility of roots in $\Pge$ is modeled by diagonals of a regular convex $(n+3)$-gon.
The negative simple roots $-\Pi=\set{-\alpha_1,\ldots,-\alpha_n}$ are identified with $n$ diagonals forming a ``snake'' configuration as illustrated in Figure~\ref{Asnake} for the case $n=7$.
The set of positive roots in type $A_n$ consists of all vectors of the form $\alpha_i+\alpha_{i+1}\cdots+\alpha_j$ for $1\le i\le j\le n$.
Each positive root $\alpha_i+\alpha_{i+1}+\cdots+\alpha_j$ is identified with the unique diagonal that crosses the diagonals labeled $-\alpha_i,-\alpha_{i+1},\ldots,-\alpha_j$ and crosses no other diagonals in the snake.

This correspondence between almost positive roots and diagonals has the property that two roots are compatible if and only if their diagonals do not cross.
Thus a cluster of almost positive roots corresponds to a choice of $n$ noncrossing diagonals, or in other words, a triangulation of the polygon.
The map $\tau_-$ acting on $\Pge$ corresponds to a symmetry of the polygon.
Specifically, $\tau_-$ acts as the unique reflective symmetry of the polygon that fixes the diagonal labeled $-\alpha_1$.
Similarly, $\tau_+$ acts as the unique reflective symmetry fixing $-\alpha_2$.
Thus $\br{\tau_+,\tau_-}$ is the full symmetry group of the polygon.

\subsection{The usual model for compatibility in type $B_n$ ($C_n$)}\label{usual BC}
Let $W$ be a Coxeter group of type $B_n$ with simple generators $S$ and bipartition of $S$ as described in Section~\ref{nc B}.
For each $i\in[n]$, let $\alpha_i$ be the simple root corresponding to $s_i$.

\begin{figure}
\begin{center} 
\scalebox{.95}{
\begin{picture}(0,0)(0,0) 
\put(26,60){$-\alpha_1$}
\put(36,130){$-\alpha_2$}
\put(78,60){$-\alpha_3$}
\put(95,130){$-\alpha_4$}
\put(140,65){$-\alpha_3$}
\put(155,140){$-\alpha_2$}
\put(187,95){--$\alpha_1$}
\end{picture} 
\includegraphics{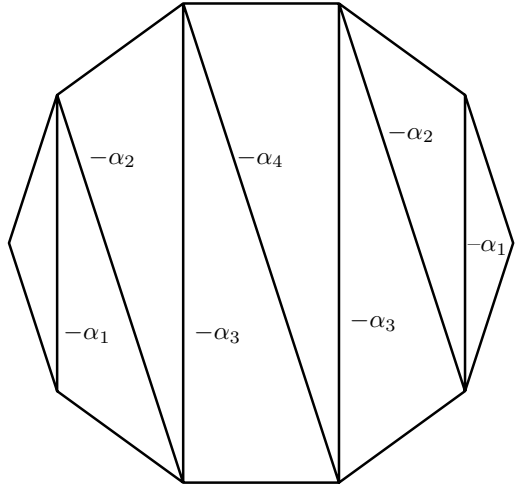}
}
\end{center} 
\caption{The ``snake'' in type $B_4$} 
\label{Bsnake} 
\end{figure}

Clusters of almost positive roots in type $B_n$ are modeled by centrally symmetric triangulations of a regular convex $(2n+2)$-gon.
Each almost positive root is identified with either a diameter or a centrally symmetric pair of diagonals and, as in type $A_n$, ``compatible'' means ``noncrossing.''
As illustrated in Figure~\ref{Bsnake}, the negative simple root $-\alpha_n$ is identified with a diameter and each other negative simple root is identified with a centrally symmetric pair of diagonals in the snake pattern.
For each negative simple root $-\alpha_s$, fix one of the (1 or 2) diagonals $\overline{QR}$ associated to $-\alpha_s$.
Then each positive root $\beta$ is assigned to the unique diameter or symmetric pair of diagonals such that the coefficient of $\alpha_s$ in the simple root expansion of $\beta$ equals the number of segments (0, 1 or 2) associated to $\beta$ that cross the diagonal $\overline{QR}$.
The almost positive roots of the dual root system of type $C_n$ are associated to the diagonals such that a root (of type $B_n$) and its coroot (of type $C_n$) are associated to the same diagonal.
The description of $\tau_\pm$ is just as in type~$A$.

\subsection{The usual model for compatibility in type $D_n$}\label{usual D}
\begin{figure}
\begin{center} 
\scalebox{.95}{
\begin{picture}(0,0)(0,0) 
\put(26,60){$-\alpha_1$}
\put(36,130){$-\alpha_2$}
\put(78,70){$-\alpha_3$}
\put(95,140){$-\alpha_5$}
\put(100,40){$-\alpha_4$}
\put(140,65){$-\alpha_3$}
\put(155,140){$-\alpha_2$}
\put(187,95){--$\alpha_1$}
\end{picture} 
\includegraphics{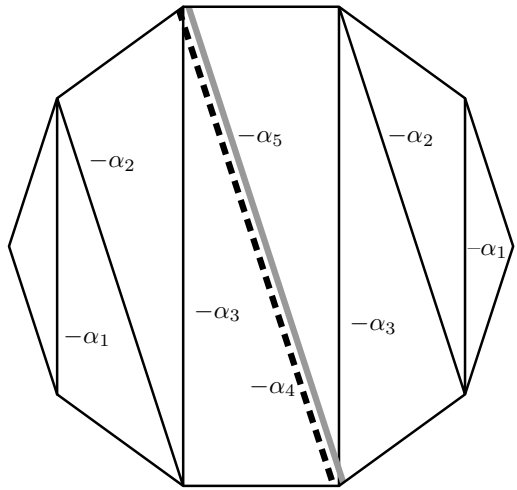}
}
\end{center} 
\caption{The ``snake'' in type $D_5$} 
\label{Dsnake} 
\end{figure}

Let $W$ be a Coxeter group of type $D_n$ with simple generators $S$ and bipartition of $S$ as described in Section~\ref{nc D}.
The model for compatibility inhabits a regular $(2n)$-gon.
Each almost positive root is identified with a colored diameter (``gray'' or ``dashed'') or with an uncolored centrally symmetric pair of diagonals.
As illustrated in Figure~\ref{Dsnake} in the case $n=5$, the simple root $-\alpha_n$ is identified with a gray diameter, $-\alpha_{n-1}$ is identified with the dashed diameter in the same location, and each other negative simple root is identified with a pair of diagonals to form the familiar snake pattern.
The remaining roots are assigned to colored diameters or pairs of uncolored diagonals as described in \cite[Section~3.5]{ga}.
Compatibility between two pairs of symmetric diagonals or between a diameter and a pair of symmetric diagonals is given by the usual noncrossing condition.
Compatibility between two diameters is as follows:
Two different diameters of the same color are always compatible, while two diameters of different colors are compatible if and only if they coincide (as for example $-\alpha_4$ and $-\alpha_5$ in Figure~\ref{Dsnake}).
The maps $\tau_\pm$ act not only geometrically, but also on colors.
As before, the map $\tau_-$ acts as the reflective symmetry of the polygon fixing $-\alpha_1$ and the map $\tau_+$ acts as the reflective symmetry fixing $-\alpha_2$.
In addition $\tau_{(-1)^n}$ changes the colors of all diameters and $\tau_{(-1)^{n-1}}$ fixes the colors of all diameters.

\section{Diagrams for almost positive roots} \label{root diag}
In this section we construct a planar diagram for each almost positive root.
The main idea for the construction springs from the similarity, in the classical types, between diagrams for noncrossing partitions and diagrams for clusters.
For example, in the case where $W$ is of type $A_n$, diagrams for noncrossing partitions involve an $(n+1)$-gon, while diagrams for clusters are triangulations of an $(n+3)$-gon.
The $(n+1)$ in this example is the Coxeter number $h$.
Similarly, diagrams for noncrossing partitions in types $B$ and $D$ involve an $h$-gon, while diagrams for clusters inhabit an $(h+2)$-gon.

Let $o$ be a smallest $W$-orbit.
Inspection of the classical examples described in Section~\ref{nc sec} and of Figure~\ref{proj orb} reveals that each $\br{c_+,c_-}$-orbit in $o$ either projects to a regular $h$-gon in $P$ or to the origin of $P$.

We begin by representing each reflection $t\in T$ by the diagram of the parabolic subgroup $\set{1,t}$.
This is a set of parallel line segments connecting various points in the projected orbit.
Some of the segments are \emph{labeled} in the following sense:
Each line segment with an endpoint at the origin is the projection of a specific line segment $L$ in $V$, and the projected segment in $P$ is considered to be labeled by the endpoint of $L$ that projects to the origin.
In particular, identical line segments in $P$ arising from different line segments in the ambient space are considered to be distinguishable.
Since each reflecting hyperplane intersects $P$ in a line, in particular no reflecting hyperplane contains $P$.
Thus no segment $L$ whose endpoints are distinct points related by a reflection in $W$ can project degenerately to the origin.

\begin{figure}
\scalebox{.8}{\includegraphics{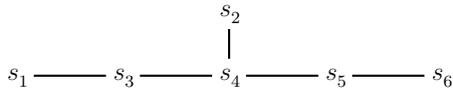}}
\caption{The Coxeter diagram of type $E_6$.}
\label{e6diag}
\end{figure}

\begin{figure}
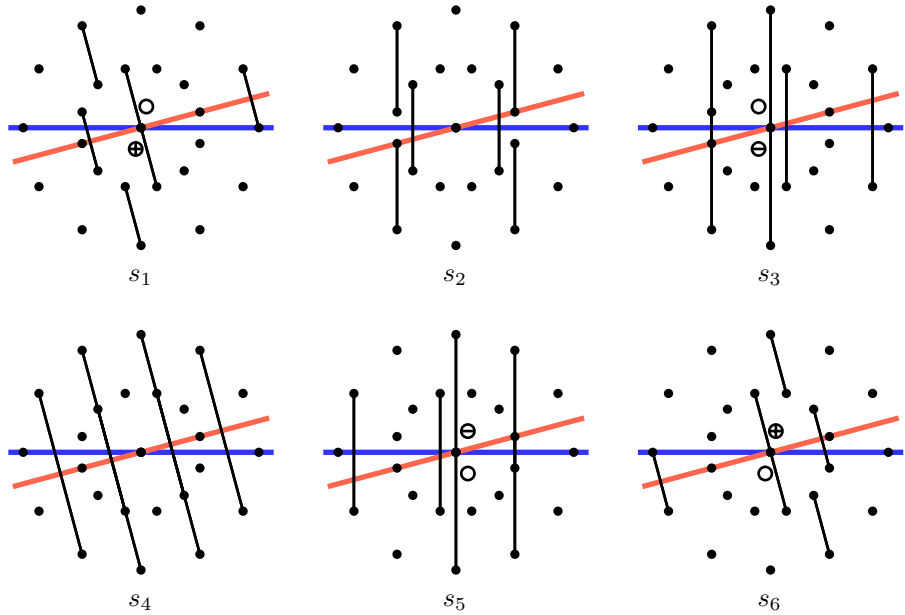

\begin{tabular}{ccc}
\figpage{ps: 1}\scalebox{1}{\includegraphics{e6simples.ps}}&
\figpage{ps: 2}\scalebox{1}{\includegraphics{e6simples.ps}}&
\figpage{ps: 3}\scalebox{1}{\includegraphics{e6simples.ps}}\\
$s_1$&$s_2$&$s_3$\\[15 pt]
\figpage{ps: 4}\scalebox{1}{\includegraphics{e6simples.ps}}&
\figpage{ps: 5}\scalebox{1}{\includegraphics{e6simples.ps}}&
\figpage{ps: 6}\scalebox{1}{\includegraphics{e6simples.ps}}\\
$s_4$&$s_5$&$s_6$
\end{tabular}
\caption{Collections of segments representing simple reflections in $E_6$}
\label{e6simp}
\end{figure}

We explain the general construction with a running example in $E_6$.
We take the bipartition $S_+=\set{s_1,s_4,s_6}$, $S_-=\set{s_2,s_3,s_5}$ where the diagram of $E_6$ is labeled as shown in Figure~\ref{e6diag}.
Figure~\ref{proj orb} shows the projection of a smallest $E_6$-orbit $o$ into the Coxeter plane, including three points of which project to the origin.
We refer to these these points as $\oplus$, \raisebox{.7 pt}{\scalebox{.72}{$\bigcirc$}}, and $\ominus$.
Figure~\ref{e6simp} shows the collections of segments, labeled as necessary to distinguish points at the origin, representing each simple reflection.

\begin{figure}
\begin{tabular}{ccc}
\figpage{ps: 2}\scalebox{1}{\includegraphics{e6clus_con.ps}}&
\figpage{ps: 3}\scalebox{1}{\includegraphics{e6clus_con.ps}}&
\scalebox{1}{\includegraphics{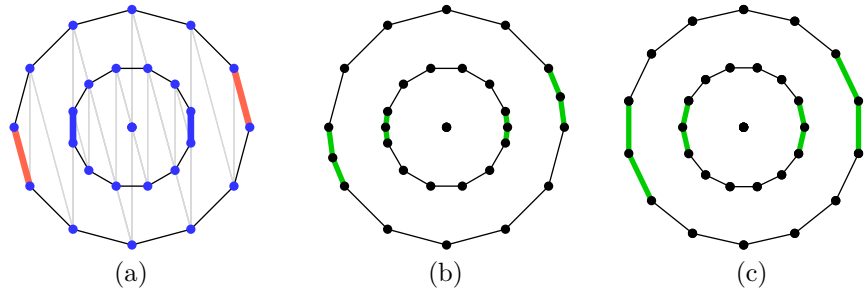}}\\
(a)&(b)&(c)
\end{tabular}
\caption{Making $(h+2)$-gons from $h$-gons in $E_6$}
\label{e6clus_con}
\end{figure}

The line segments representing simple reflections $-\alpha_i$ for $s_i\in S_+$ are perpendicular to the line $L_+$ and the line segments representing simple reflections $-\alpha_i$ for $s_i\in S_-$ are perpendicular to the line $L_-$.
(These are the lines $L_\pm$ from Section~\ref{par sec}.)
In Figure~\ref{e6simp}, the red (off-horizontal) line is $L_+$ and the blue (horizontal) line is $L_-$.
Reflecting a diagram through the line $L_\ep$ corresponds to conjugating a reflection (or a parabolic subgroup) by $c_\ep$.
Call an edge on an $h$-gon in the Coxeter plane \emph{distinguished} if it is perpendicular to $L_+$ or $L_-$.
The distinguished edges are shown in red and blue in Figure~\ref{e6clus_con}.a.

The $h$-gons in the diagram are then expanded into $(h+2)$-gons as follows:  for each distinguished edge, a new vertex is placed in convex position with the other vertices, in such a way that the new vertex has the effect of subdividing the distinguished edge, as shown in Figure~\ref{e6clus_con}.b, with the new edges shown in green.
Finally, each $(h+2)$-gon is made regular, while preserving the distances of points to the origin.
Also preserved are the relative orientations of the polygons:  any two of the $h$-gons either interlace, like the two $h$-gons for $E_6$, or line up, like any pair of $h$-gons for $E_7$.
 (See Figure~\ref{proj orb}.)
The corresponding $(h+2)$-gons are made to line up or interlace accordingly.
This construction of one regular $(h+2)$-gon for each regular $h$-gon in the projected orbit is unique up to a global scaling and rotation.
Any multiple points at the origin in the projected orbit remain at the origin with the same multiplicity.
Figure~\ref{e6clus_con}.c shows the regular $(h+2)$-gons in the $E_6$ example.

\begin{figure}
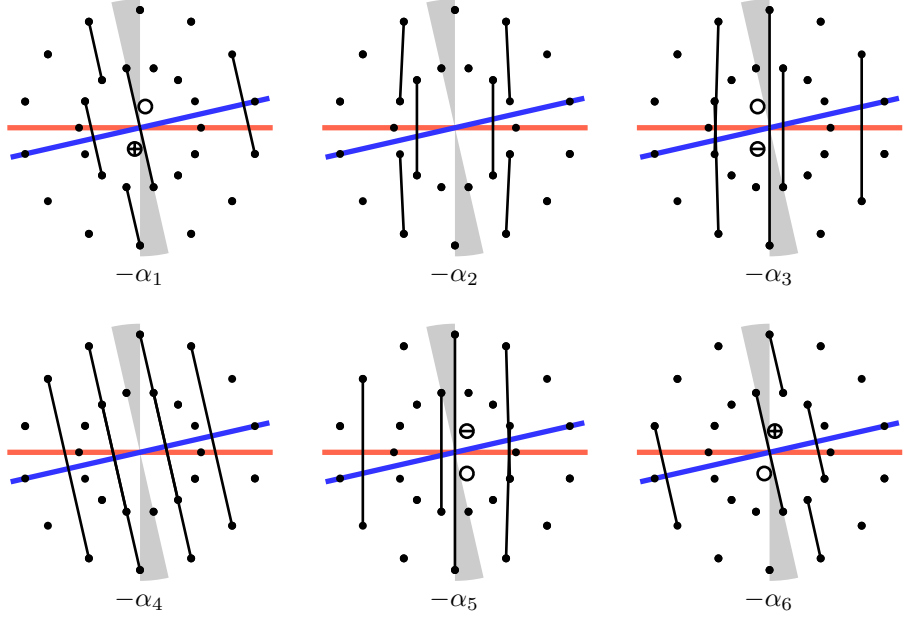

\begin{tabular}{ccc}
\figpage{ps: 1}\scalebox{1}{\includegraphics{e6negsimples.ps}}&
\figpage{ps: 2}\scalebox{1}{\includegraphics{e6negsimples.ps}}&
\figpage{ps: 3}\scalebox{1}{\includegraphics{e6negsimples.ps}}\\
$-\alpha_1$&$-\alpha_2$&$-\alpha_3$\\[15 pt]
\figpage{ps: 4}\scalebox{1}{\includegraphics{e6negsimples.ps}}&
\figpage{ps: 5}\scalebox{1}{\includegraphics{e6negsimples.ps}}&
\figpage{ps: 6}\scalebox{1}{\includegraphics{e6negsimples.ps}}\\
$-\alpha_4$&$-\alpha_5$&$-\alpha_6$
\end{tabular}
\caption{Collections of segments representing negative simple roots in $E_6$}
\label{e6negsimp}
\end{figure}

Once regular $(h+2)$-gons are constructed, the collections of line segments representing simple reflections in the original projection are used to represent negative simple roots.
Specifically, each simple reflection $s$ is represented by a collection of (labeled) line segments connecting pairs of vertices in the original $h$-gons.
Our construction defines an inclusion from the (labeled) vertices of the projected orbit to the points of the new diagram.
We use this inclusion to define a collection of line segments representing the negative simple root $-\alpha_s$.
Figure~\ref{e6negsimp} shows the collections of labeled line segments representing the negative simple roots of $E_6$.
Notice that some collinearities of points present in the diagrams for simple reflections are not present in the diagrams for negative simple roots.

In the case of planar diagrams for partitions, projection to the Coxeter plane defines a collection of line segments representing every reflection (and more generally every parabolic subgroup).
The dihedral action on the diagrams given by symmetries of the $h$-gons corresponds to the conjugation action of $\br{c_+,c_-}$ on reflections.
In contrast, in our construction of planar diagrams for almost positive roots, once the projected orbit has been altered as described above, there is no longer an identification, via projections, of positive roots with collections of segments.
To map each positive root to a collection of segments, we first need an appropriate action of $\br{\tau_+,\tau_-}$ on labeled line segments.
The action is the dihedral action given by symmetries of the $(h+2)$-gons, with alterations involving multiple points at the origin.

Let $L'_+$ be the line through the origin such that reflection in $L'_+$ preserves the diagrams representing negative simple roots $-\alpha_i$ for $s_i\in S_-$.
Let $L'_-$ be the line such that reflection in $L'_-$ preserves the diagrams representing negative simple roots $-\alpha_i$ for $s_i\in S_+$.
(Note the mixture of signs in the definitions of $L'_+$ and $L_-$.)
Figure~\ref{e6negsimp} shows $L'_+$ in red and $L'_-$ in blue.

Suppose $AB$ is a line segment in a diagram for an almost positive root.
If neither $A$ nor $B$ lies at the origin then $\tau_\ep$ acts on the segment by acting on each endpoint as a reflection through $L'_\ep$.
If one of the endpoints (say $A$) lies at the origin then, in most cases, $\tau_\ep$ acts on $B$ as a reflection through $L'_\ep$ and acts on $A$ by the action of $c_\ep$, which permutes the origin points.
The exception is that, when $A$ is at the origin and $AB$ is perpendicular to $L'_\ep$, the segment $AB$ is fixed by $\tau_\ep$.
This exceptional case occurs only when $AB$ is one of the segments representing some $-\alpha_i$ for $s_i\in S_{-\ep}$.
In particular, $\tau_\ep$ preserves the entire diagram of each $-\alpha_i$ for $s_i\in S_{-\ep}$.

\begin{figure}
\includegraphics{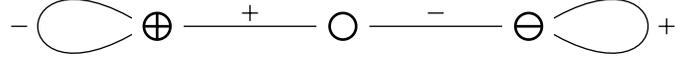}
\caption{The action of $\br{c_+,c_-}$ on origin points in $E_6$}
\label{three}
\end{figure}

\begin{figure}
\centerline{\scalebox{.959}{\includegraphics{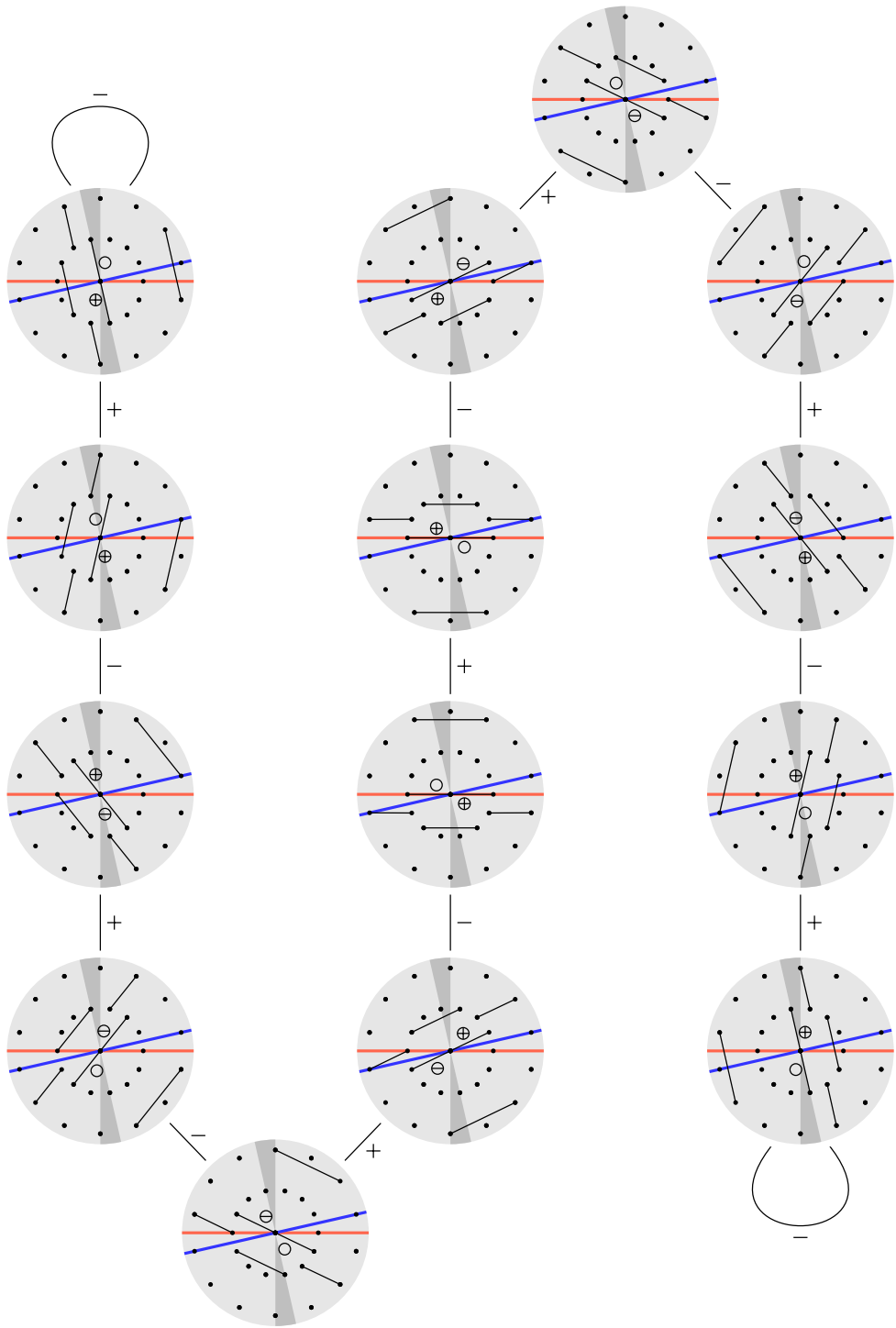}}}
\caption{The $\br{\tau_+,\tau_-}$-orbit of $-\alpha_1$ (top left) and $-\alpha_6$ (bottom right) in type $E_6$}
\label{e6tau}
\end{figure}

In the $E_6$ example, the group $\br{c_+,c_-}$ acts on the three points $\oplus$, \raisebox{.7 pt}{\scalebox{.72}{$\bigcirc$}}, and $\ominus$ at the origin according to the diagram in Figure~\ref{three}, where edges marked $+$ or $-$ denote the action of $c_+$ or $c_-$ respectively.
Figure~\ref{e6tau} shows, in type $E_6$, the $\br{\tau_+,\tau_-}$-orbit of the diagram of the negative simple roots $-\alpha_1$ and $-\alpha_6$.

Having defined $\tau_+$ and $\tau_-$, we require for any almost positive root $\beta$ and $\ep\in\set{+.-}$ that the diagram for $\beta$ and the diagram for $\tau_\ep(\beta)$ be related by the action of $\tau_\ep$ on diagrams.
In light of Propositions~\ref{c orbits} and~\ref{tau orbits}, this requirement uniquely defines a collection of segments representing each almost positive root.

The assignment of diagrams to almost positive roots can also be easily explained in terms of the transformation $\tau_+\tau_-$.
The action of $\tau_+\tau_-$ on $P$ is a rotation of order $h+2$.
To describe the action of $\tau_+\tau_-$ on the labels of points at the origin, we define the \emph{gray zone} to be the union of the two acute wedges defined by two lines through the origin:  one perpendicular to $L'_+$ and the other perpendicular to $L'_-$.
The gray zone is shown in the almost positive root diagrams of Figures~\ref{e6negsimp} and~\ref{e6tau}.
When $\tau_+\tau_-$ acts on a segment containing an origin point, the non-origin point is rotated while $c_+c_-$ is applied to the label of the point at the origin.
The exception to this rule occurs when the segment is rotated from one side of the gray zone to the other:  in this case the label of the origin point is fixed.
Inspection of Figure~\ref{e6tau} reveals why the label of the origin point is fixed when the segment is rotated past the gray zone.
(In the figure, this rule applies when $\tau_+\tau_-$ is applied to $-\alpha_1$ or to $-\alpha_6$.)
When a segment contains an origin point and is incident to the gray zone, one of $\tau_\pm$ (in this case $\tau_+$) fixes the segment, while the other (in this case $\tau_-$) reflects the segment while swapping the label of the origin point with the label of the origin point of the antipodal segment.
The swap of labels of the origin points must occur because the conjugation action of $c_\ep$ fixes a simple reflection $\alpha_i$ when $s_i\in S_\ep$.

\section{Criteria for compatibility}\label{cl diag}
We now give a type-by-type account of the construction of Section~\ref{root diag}.
Computer data for each type is available as described in Section~\ref{computer}.
In each type, we attempt to define a criterion for compatibility of almost positive roots.
In types $A_n$, $B_n$, and $I_2(m)$, the criterion is a simple noncrossing condition.
Comparing the type-$D$ construction with the model described in Section~\ref{usual D}, we obtain a simple rule for compatibility in the presence of multiple origin points.
Inspection of the group of type $F_4$ reveals a criterion for compatibility in the presence of multiple $(h+2)$-gons, and, surprisingly, this criterion carries forward to the much more complicated diagrams of type $H_4$.
By combining the multiple-origin-points rule from type $D$ and the multiple-polygons rule from $F_4$, one hopes to obtain a criterion that works in generality.
Indeed, the criterion works in type $E_6$, but breaks down when applied to type $E_7$.
The criterion also fails for $E_8$, whose diagrams have no origin points.

In each of the types $E_6$, $E_7$, $E_8$, $F_4$, and $H_4$, we make \textit{ad hoc} alterations of the Coxeter-plane diagrams which yield a model with a simple noncrossing condition on compatibility.
In each case, alterations are made to the diagrams for negative simple roots, and then carried over to the positive roots by the action of $\br{\tau_+,\tau_-}$.

\subsection{Type $A_n$}\label{cl A}
\begin{figure}
\centerline{\scalebox{.9}{
\begin{picture}(0,0)(0,0) 
\put(4,65){$s_1$}
\put(28,65){$s_2$}
\put(54,60){$s_3$}
\put(81,70){$s_4$}
\put(108,65){$s_5$}
\put(122,98){$s_6$}
\put(145,65){$s_7$}
\end{picture} 
\includegraphics{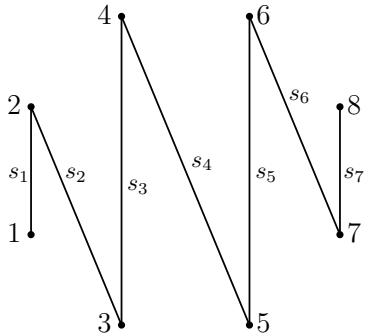}}
}
\caption{The action of simple reflections on the the smallest orbit of $A_7$}
\label{Acsnake}
\end{figure}

We continue the notation of Section~\ref{nc A}.
The orbit $\set{e_i:i\in[n+1]}$ projects into $P$ as a regular $(n+1)$-gon with vertices labeled according to the $(n+1)$-cycle $c$.
Each reflection in $W$ fixes all but two elements of the orbit, which are transposed.
The simple reflections $s_i=(i\,\,i\!+\!1)$ act on the orbit as illustrated in Figure~\ref{Acsnake} for $n=7$, and thus the ``snake'' configuration (see Figure~\ref{Asnake}) of Fomin and Zelevinsky's model arises from the construction described in Section~\ref{root diag}, rather than as an \textit{ad hoc} construction.
In the usual type-$A$ model, the action of $\br{\tau_+,\tau_-}$ on almost positive roots corresponds to the dihedral action on $(n+3)$-gons.
In the Coxeter-plane model, the same correspondence holds by construction.
Thus the criterion for compatibility of diagrams in the Coxeter plane is the usual model:
\begin{cl_criterion}\label{cl simple}
Two distinct almost positive roots, represented by collections $C_1$ and $C_2$ of line segments in $P$, are compatible if and only if no segment in $C_1$ crosses any segment in $C_2$.
\end{cl_criterion}

\subsection{Type $B_n$}\label{cl B}
Arguing as in Section~\ref{cl A}, the type-$B$ case also reduces to the model from~\cite{ga}, described in Section~\ref{usual BC}, so Compatibility Criterion~\ref{cl simple} describes compatibility.

\subsection{Type $I_2(m)$}\label{cl I}
\begin{figure}
\centerline{\scalebox{.75}{\includegraphics{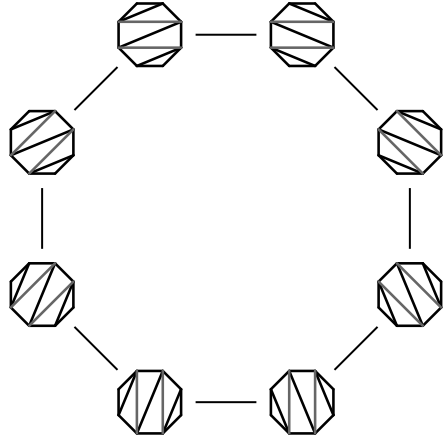}}}
\caption{The $G_2$ associahedron.}
\label{g2assoc}
\end{figure}

In this case, $h=m$ and~$W$ has an orbit of size $m$.
The projected orbit is an $m$-gon, which is deformed to an $(m+2)$-gon.
There are $m+2$ almost positive roots, and each corresponds to a maximal set of parallel diagonals of the $(m+2)$-gon.
One easily checks that the criterion for compatibility is again Compatibility Criterion~\ref{cl simple}.
The clusters are shown in Figure~\ref{g2assoc} for the case $I_2(6)=G_2$.
The edges shown connecting clusters are the edges of the corresponding associahedron.

\subsection{Type $H_3$}\label{cl H3}
\begin{figure}
\psfrag{1}[cc][cc]{$-\alpha_1$}
\psfrag{2}[cc][cc]{$-\alpha_2$}
\psfrag{3}[cc][cc]{$-\alpha_3$}
\centerline{\scalebox{.9}{\includegraphics{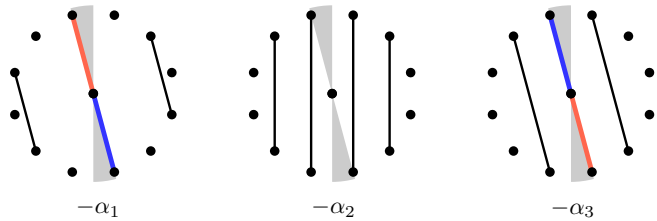}}}
\caption{Negative simple roots in $H_3$.}
\label{h3roots}
\end{figure}

As discussed in Section~\ref{nc H3}, a Coxeter group of type $H_3$ has $h=10$ and has smallest orbit of size $12$.
Coxeter-plane diagrams for the negative simple roots are shown in Figure~\ref{h3roots}.
(Cf.\ the bottom row of Figure~\ref{h3nc}.)
Two of the negative simple roots include segments connecting to origin points. 
To determine a criterion for compatibility, we check all $3\cdot 17=51$ pairs $-\alpha_s,\beta$ where $\alpha_s$ is a negative simple root and $\beta$ is an almost-positive root not equal to $-\alpha_s$.
The following simple criterion is correct:

\begin{cl_criterion}\label{cl coincide}
Two distinct almost positive roots, represented by collections $C_1$ and $C_2$ of line segments in $P$, are compatible if and only if no segment in $C_1$ crosses any segment in $C_2$, unless the two segments coincide (as unlabeled segments in $P$).
\end{cl_criterion}
Compatibility Criterion~\ref{cl coincide} also defines compatibility in the cases $A_n$, $B_n$, and $I_2(m)$.
In these cases, there are no coinciding segments in the representations of distinct roots because there are no multiple points in the projected orbit.

\begin{figure}
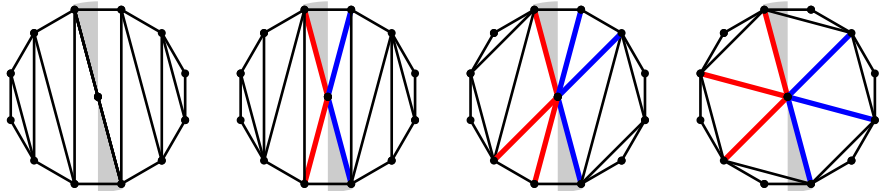

\begin{tabular}{cccc}
\figpage{ps: 4}\includegraphics{h3clusters.ps}&
\figpage{ps: 3}\includegraphics{h3clusters.ps}&
\figpage{ps: 2}\includegraphics{h3clusters.ps}&
\figpage{ps: 1}\includegraphics{h3clusters.ps}
\end{tabular}
\caption{Diagrams for clusters of type $H_3$.}
\label{h3clus}
\end{figure}

A representative of each $\br{\tau_+,\tau_-}$-orbit of clusters is shown in Figure~\ref{h3clus}.
In this figure, the edges of the polygon are also drawn.
The labeling of the origin points is irrelevant to Compatibility Criterion~\ref{cl coincide}, but is indicated in Figure~\ref{h3clus} by colors as in Figure~\ref{h3roots}.
The leftmost diagram in Figure~\ref{h3clus} includes the roots $-\alpha_1$ and $-\alpha_3$, and thus has two red-blue pairs of coinciding segments attached to the origin.  
These pairs of segments are shown in black.

\subsection{Type $D_n$}
\begin{figure}
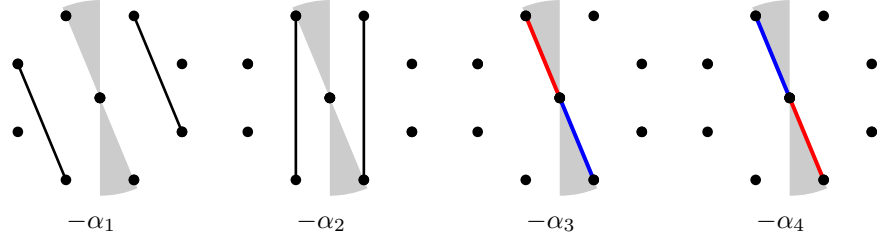

\begin{tabular}{cccc}
\figpage{ps: 1}\includegraphics{d4negsimples.ps}&
\figpage{ps: 2}\includegraphics{d4negsimples.ps}&
\figpage{ps: 3}\includegraphics{d4negsimples.ps}&
\figpage{ps: 4}\includegraphics{d4negsimples.ps}\\
$-\alpha_1$&$-\alpha_2$&$-\alpha_3$&$-\alpha_4$
\end{tabular}
\caption{Diagrams for almost positive roots of type $D_4$.}
\label{d4negsimp}
\end{figure}

We continue the notation and conventions of Section~\ref{nc D}.
Diagrams for the negative simple roots of type $D_4$, shown in Figure~\ref{d4negsimp}, are representative of the diagrams for general $n$.

To compare the Coxeter-plane diagrams for almost positive roots to the usual model described in Section~\ref{usual D}, we color the origin points red and blue so that the diagram for $-\alpha_{n-1}$ has a red segment lying to the left of the gray zone and a blue segment lying to the right of the gray zone.
The diagram for $-\alpha_n$ will then have a blue segment lying to the left of the gray zone and a red segment lying to the right of the gray zone.
The Coxeter-plane diagrams for almost positive roots correspond to the usual diagrams described in Section~\ref{usual D}, equating a pair of half diameters with a (full) dashed diameter if the red half-diameter is to the left of the gray zone, and to a (full) gray diameter if the blue half-diameter is to the right of the gray zone.

The origin points were colored red and blue so as to make the correspondence hold for negative simple roots.
Thus, to verify the correspondence, it remains only to show that the action of $\tau_+$ and $\tau_-$ preserves the correspondence.
The case of non-diameters is trivial.
The two origin points arise as the images of $\pm e_n$ under projection to $P$.
If $n$ is even, then $c_+=s_1s_3\cdots s_{n-1}\cdot s_n$, which transposes $e_n$ and $-e_n$, while $c_-=s_2s_4\cdots s_{n-2}$, which fixes both $e_n$ and $-e_n$.
Since $s_{n-1}$ and $s_n$ are both in $S_+$, the exception in the definition of the action of $\tau_\ep$ does not apply to $\tau_+$.
Thus $\tau_+$ reflects a half-diameter, keeping it on the same side of the gray zone, and changing its origin point.
This corresponds to the action of $\tau_+$ on colored diameters, which applies the same reflection and reverses the color of the diameter.
Similarly, $\tau_-$ reflects a half-diameter, keeping it on the same side of the gray zone, and fixes its origin point, corresponding to the action of $\tau_-$ on colored diameters, which applies the reflection and fixes the color of diameters.
If $n$ is odd, the correspondence works out similarly.
The compatibility criterion in the usual $D_n$ model translates to the following criterion for compatibility for Coxeter-plane diagrams,
which is also correct in types $A_n$, $B_n$, $I_2(m)$, and~$H_3$.

\begin{cl_criterion}\label{cl coincide center}
Two distinct almost positive roots, represented by collections $C_1$ and $C_2$ of line segments in $P$, are compatible if and only if both of the following conditions holds:
\begin{enumerate}
\item No segment in $C_1$ crosses any segment in $C_2$, unless the two segments coincide (as unlabeled segments in $P$).
\item If two segments---one from $C_1$ and one from $C_2$---involve the same origin point then either both segments are contained in a common line or the union of the two segments does not cross the gray zone.  
\end{enumerate}
\end{cl_criterion}

\subsection{Type $F_4$}
\begin{figure}
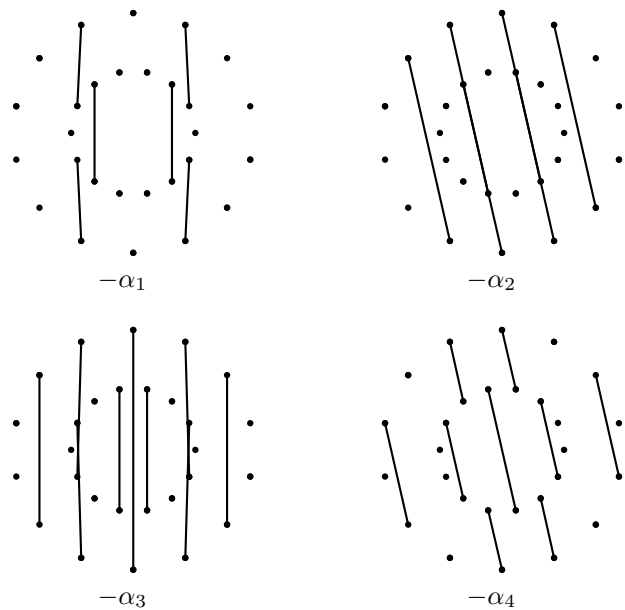

\begin{tabular}{ccccc}
\scalebox{.8}{\figpage{ps: 1}\includegraphics{f4root.ps}}&&&&\scalebox{.8}{\figpage{ps: 2}\includegraphics{f4root.ps}}\\[-2 pt]
$-\alpha_1$&&&&$-\alpha_2$\\[10 pt]
\scalebox{.8}{\figpage{ps: 3}\includegraphics{f4root.ps}}&&&&\scalebox{.8}{\figpage{ps: 4}\includegraphics{f4root.ps}}\\[-2 pt]
$-\alpha_3$&&&&$-\alpha_4$
\end{tabular}
\caption{Diagrams for negative simple roots in $F_4$}
\label{f4roots}
\end{figure}

The Coxeter group of type $F_4$ has Coxeter number $h=12$ and smallest orbit of size $24$.
We label the diagram in a linear order and set $S_+=\set{s_2,s_4}$ and $S_-=\set{s_1,s_3}$.
The projected orbit consists of two $12$-gons, which become two $14$-gons in the cluster diagrams.
Figure~\ref{f4roots} shows diagrams for the negative simple roots in $F_4$.
The diagram for $-\alpha_2$ has two sets of four nearly collinear points.
Ordering the points from top to bottom, the diagram contains a segment connecting the first and third points and a segment connecting the second and fourth points.
The diagram for $-\alpha_3$ also has two sets of four nearly collinear points, connected by segments in the same way.

\begin{figure}[p]
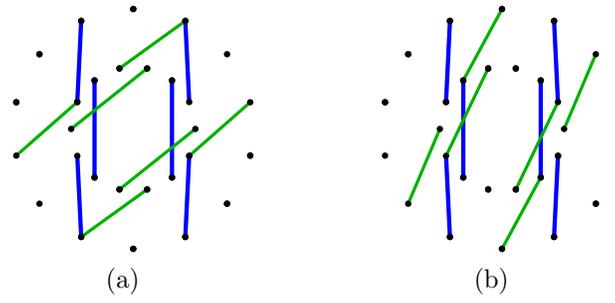

\begin{tabular}{ccccc}
\scalebox{.8}{\figpage{ps: 1}\includegraphics{f4pair.ps}}&&&&\scalebox{.8}{\figpage{ps: 2}\includegraphics{f4pair.ps}}\\
(a)&&&&(b)
\end{tabular}
\caption{a:  Coxeter-plane diagram for two compatible roots in~$F_4$. 
b:  Coxeter-plane diagram for two incompatible roots in $F_4$.}
\label{f4 compat incompat}
\end{figure}

It is immediately apparent that some crossing of segments must be allowed in the criterion for compatibility.
For example, the negative simple roots $-\alpha_3$ and $-\alpha_4$ (pictured in Figure~\ref{f4roots}) are compatible.
Figure~\ref{f4 compat incompat} illustrates a crossing that must be allowed, and a crossing that must be ruled out.
\begin{figure}[p]
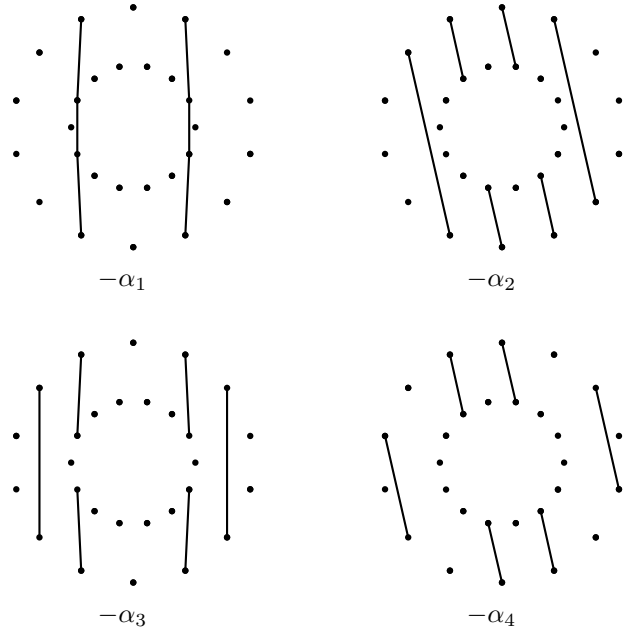

\begin{tabular}{ccccc}
\scalebox{.8}{\figpage{ps: 5}\includegraphics{f4root.ps}}&&&&\scalebox{.8}{\figpage{ps: 6}\includegraphics{f4root.ps}}\\
$-\alpha_1$&&&&$-\alpha_2$\\[15 pt]
\scalebox{.8}{\figpage{ps: 7}\includegraphics{f4root.ps}}&&&&\scalebox{.8}{\figpage{ps: 8}\includegraphics{f4root.ps}}\\
$-\alpha_3$&&&&$-\alpha_4$
\end{tabular}
\caption{Altered negative simple roots in $F_4$}
\label{f4roots_alt}
\end{figure}
\begin{figure}[p]
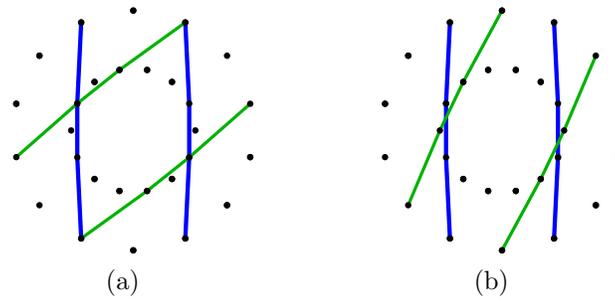

\begin{tabular}{ccccc}
\scalebox{.8}{\figpage{ps: 3}\includegraphics{f4pair.ps}}&&&&\scalebox{.8}{\figpage{ps: 4}\includegraphics{f4pair.ps}}\\
(a)&&&&(b)
\end{tabular}
\caption{a:  Altered diagram for two compatible roots in~$F_4$. 
b:  Altered diagram for two incompatible roots in $F_4$.}
\label{f4 compat incompat alt}
\end{figure}
In each picture, the blue segments represent one of the roots and the green segments represent the other.

Compatibility appears to depend not only whether two segments cross, but where they cross, and what other segments they are incident to.
Specifically, given a segment connecting points in the cluster diagram, the \emph{outer ring} of the segment is the larger of the two $(h+2)$-gons containing the endpoints of the segment. 
(If the two endpoints are on the same $(h+2)$-gon then that $(h+2)$-gon is the outer ring.)
Given a pair of segments, the outer ring of the pair is the larger of the two outer rings.
Given the diagrams for two positive roots, a segment is \emph{active} with respect to the two roots if its outer ring $R$ is the outermost $(h+2)$-gon of the diagram or if it shares a point with an active segment whose outer ring is the next larger ring than $R$.
Compatibility in $F_4$ can be described by the following criterion:

\begin{cl_criterion}\label{cl rings}
Two distinct almost positive roots, represented by collections $C_1$ and $C_2$ of line segments in $P$, are compatible if and only if no two active segments cross within the annulus between the outer ring of the pair and the next smaller $(h+2)$-gon.
\end{cl_criterion}

Compatibility Criterion~\ref{cl rings} declares the two roots in Figure~\ref{f4 compat incompat}.a to be compatible because the two segments that cross are not active.
The two roots in Figure~\ref{f4 compat incompat}.b are declared to be incompatible because two active segments cross.

The cluster diagrams for $F_4$ can be altered so as to give a simpler criterion.
The alteration has the disadvantage of being \textit{ad hoc}, but has the advantage of producing simpler pictures.
The altered negative simple roots are shown in Figure~\ref{f4roots_alt}.

Compatibility of almost positive roots of type $F_4$ is described by applying Compatibility Criterion~\ref{cl coincide} to the altered diagrams.
To apply the criterion correctly, it is important to notice that the diagram for $-\alpha_1$ consists of six segments, arranged in two polygonal paths.
Figure~\ref{f4 compat incompat alt} shows the altered diagrams for the two pairs of roots whose unaltered diagrams appear in Figure~\ref{f4 compat incompat}.

\subsection{Type $H_4$}
\begin{figure}
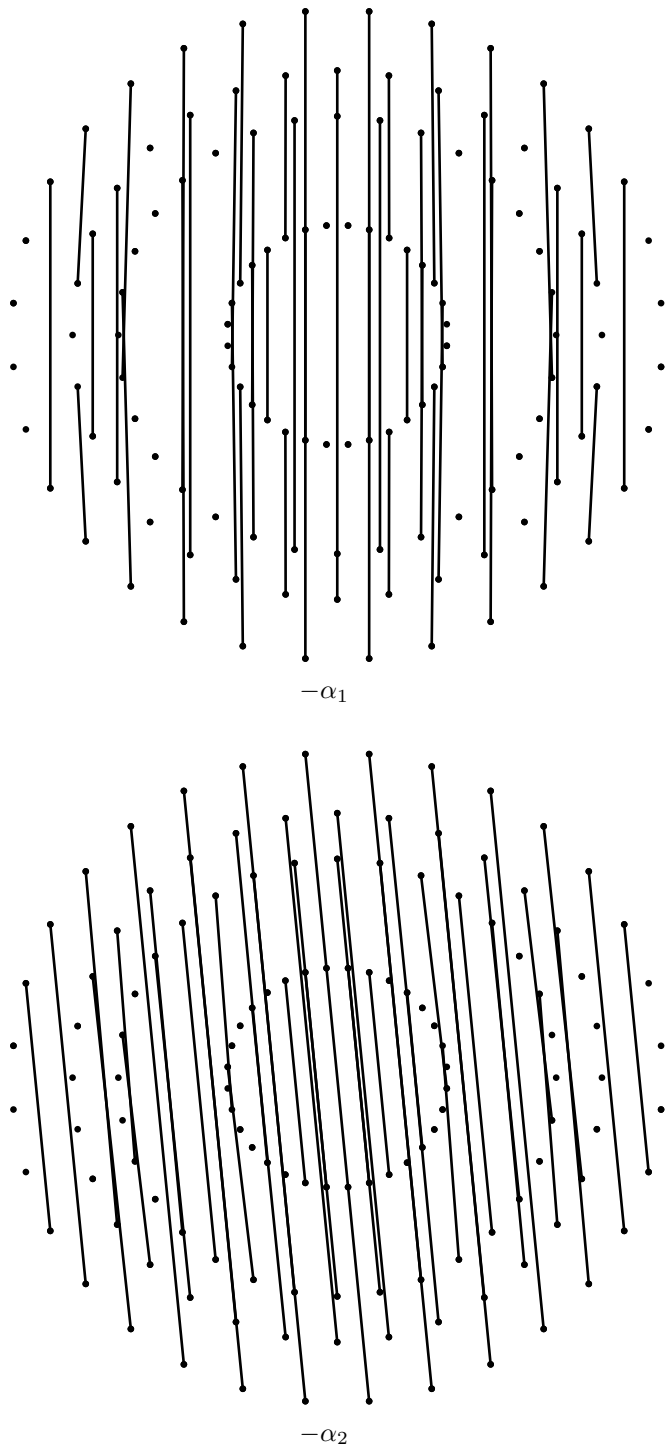

\begin{tabular}{c}
\scalebox{1}{\figpage{ps: 1}\includegraphics{h4root.ps}}\\
$-\alpha_1$\\[14.5 pt]
\scalebox{1}{\figpage{ps: 2}\includegraphics{h4root.ps}}\\
$-\alpha_2$\\[14.5 pt]
\end{tabular}
\caption{Diagrams for negative simple roots in $H_4$}
\label{h4roots}
\end{figure}

The Coxeter group of type $H_4$ has Coxeter number $h=30$ and smallest orbit of size $120$.
The projected orbit consists of four $h$-gons, which become four $32$-gons in the cluster diagrams.
Figures~\ref{h4roots} and~\ref{moreh4roots} show diagrams for the negative simple roots in $H_4$.
\begin{figure}
\begin{tabular}{c}
\scalebox{1}{\figpage{ps: 3}\includegraphics{h4root.ps}}\\
$-\alpha_3$\\[14.5 pt]
\scalebox{1}{\figpage{ps: 4}\includegraphics{h4root.ps}}\\
$-\alpha_4$\\[14.5 pt]
\end{tabular}
\caption{More diagrams for negative simple roots in $H_4$}
\label{moreh4roots}
\end{figure}

These diagrams appear almost hopelessly complex.
But surprisingly, Compatibility Criterion~\ref{cl rings} describes compatibility exactly in type $H_4$ as well.
In fact, something even nicer happens.
For every pair of incompatible roots, there is a pair of segments, each with an endpoint on the outermost $(h+2)$-gon, which cross outside the second-outermost $(h+2)$-gon.
One can alter the diagrams by erasing the parts of the segments that are inside the second-outermost $(h+2)$-gon, and then moving the remaining parts slightly so that each remaining segment connects two vertices contained in the outermost two $(h+2)$-gons.
The altered diagrams for negative simple roots are illustrated in Figure~\ref{h4rootsalt}.
Compatibility of the resulting diagrams for almost positive roots is described by Compatibility Criterion~\ref{cl coincide}.

\begin{figure}
\begin{tabular}{ccccc}
\scalebox{.5}{\figpage{ps: 5}\includegraphics{h4root.ps}}&&&&\scalebox{.5}{\figpage{ps: 6}\includegraphics{h4root.ps}}\\
$-\alpha_1$&&&&$-\alpha_2$\\[15 pt]
\scalebox{.5}{\figpage{ps: 7}\includegraphics{h4root.ps}}&&&&\scalebox{.5}{\figpage{ps: 8}\includegraphics{h4root.ps}}\\
$-\alpha_3$&&&&$-\alpha_4$
\end{tabular}
\caption{Altered diagrams for negative simple roots in $H_4$}
\label{h4rootsalt}
\end{figure}

\subsection{Type $E_6$}
The construction of Coxeter-plane diagrams for almost positive roots of type $E_6$ was discussed in detail in Section~\ref{root diag}.
The diagrams constructed feature both multiple $(h+2)$-gons and multiple origin points.
The following is a natural common extension of Compatibility Criteria~\ref{cl coincide center} and~\ref{cl rings}.

\begin{cl_criterion}\label{cl coincide center active}
Let $\alpha_1,\alpha_2$ be distinct almost positive roots, represented by collections $C_1$ and $C_2$ of line segments in $P$.
Then $\alpha_1$ and $\alpha_2$ are compatible if and only if both of the following conditions holds:
\begin{enumerate}
\item No two active segments cross within the annulus between the outer ring of the pair and the next smaller $(h+2)$-gon.
\item If two active segments---one from $C_1$ and one from $C_2$---connect the same origin point to points on the innermost $(h+2)$-gon, then either both segments are contained in a common line or the union of the two segments does not cross the gray zone.  
\end{enumerate}
\end{cl_criterion}

Compatibility in type $E_6$ is indeed described by Compatibility Criterion~\ref{cl coincide center active}.
But in fact, every Coxeter-plane diagram of a pair of incompatible roots exhibits a violation of condition (1) of Compatibility Criterion~\ref{cl coincide center active}.
In other words, compatibility in $E_6$ is also described by Compatibility Criterion~\ref{cl rings}.

\begin{figure}
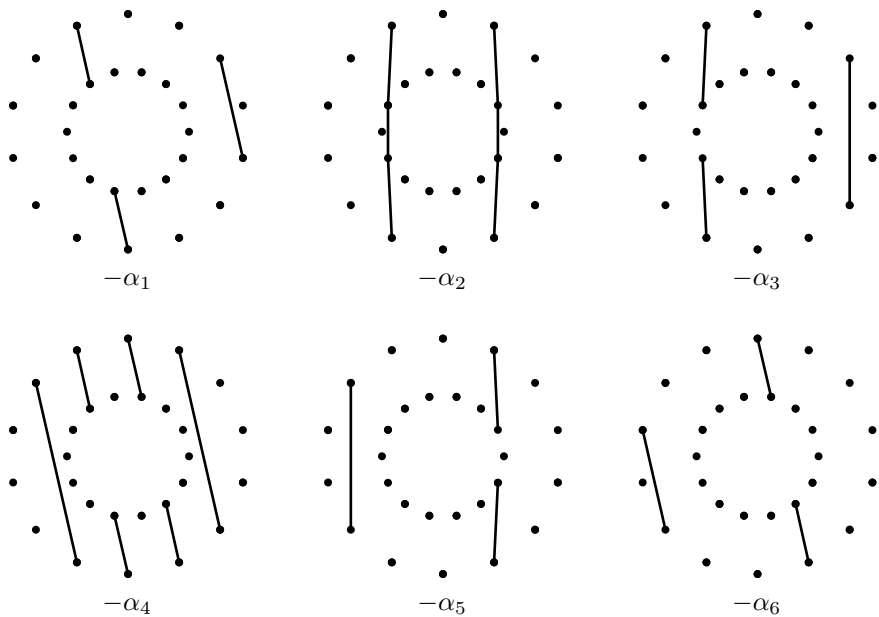

\begin{tabular}{ccc}
\figpage{ps: 7}\scalebox{1}{\includegraphics{e6negsimples.ps}}&
\figpage{ps: 8}\scalebox{1}{\includegraphics{e6negsimples.ps}}&
\figpage{ps: 9}\scalebox{1}{\includegraphics{e6negsimples.ps}}\\
$-\alpha_1$&$-\alpha_2$&$-\alpha_3$\\[15 pt]
\figpage{ps: 10}\scalebox{1}{\includegraphics{e6negsimples.ps}}&
\figpage{ps: 11}\scalebox{1}{\includegraphics{e6negsimples.ps}}&
\figpage{ps: 12}\scalebox{1}{\includegraphics{e6negsimples.ps}}\\
$-\alpha_4$&$-\alpha_5$&$-\alpha_6$
\end{tabular}
\caption{Altered diagrams for negative simple roots in $E_6$}
\label{e6negsimp alt}
\end{figure}

The $E_6$ diagrams can also be altered to produce diagrams for which compatibility is described by Compatibility Criterion~\ref{cl coincide}.
The altered simple roots are shown in Figure~\ref{e6negsimp alt}.

\subsection{Types $E_7$ and $E_8$}
\begin{figure}
\begin{tabular}{cccccc}
\scalebox{1}{\figpage{ps: 15}\includegraphics{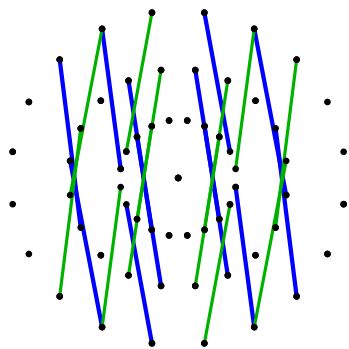}}
\end{tabular}
\caption{A counterexample to Compatibility Criterion~\ref{cl coincide center active} in $E_7$}
\label{e7counter}
\end{figure}

Unfortunately, Compatibility Criterion~\ref{cl coincide center active} breaks down in type $E_7$.
Figure~\ref{e7counter} shows the diagrams for the almost positive roots $-\alpha_3$ and $\alpha_3$, which are incompatible, but whose diagrams would be declared compatible by Compatibility Criterion~\ref{cl coincide center active}.
The counterexample seems to arise from the multiplicity of $(h+2)$-gons, rather than from the multiplicity of origin points.  
Thus it is not surprising that a similar counterexample occurs in $E_8$. 
\textit{Ad hoc} alterations of the Coxeter-plane diagrams yield diagrams for almost positive roots of types $E_7$ and $E_8$ for which compatibility is described by Compatibility Criterion~\ref{cl coincide}.

\section{Acknowledgments}
The author wishes to thank Sergey Fomin, Jon McCammond, Frank Sottile and John Stembridge for helpful conversations.

\end{document}